\documentclass[a4paper,12pt]{amsart}
\usepackage[]{amscd,amsfonts,amsmath,amsxtra,amssymb,bm}
\usepackage{graphics}
\usepackage{extarrows}
\usepackage[all]{xy}
\usepackage{hyperref}
\usepackage{eucal}
\usepackage{rotating}
\usepackage[T1]{fontenc}
\newtheorem{sub}{}[section]
\newtheorem{subsub}{}[sub]
\usepackage{rotating}
\usepackage[active]{srcltx}
\def\ov#1{\overline{#1}}

\def\coker{\mathop{\rm coker}\nolimits}
\def\Hom{\mathop{\rm Hom}\nolimits}

\def\HHom{\mathop{\mathcal Hom}\nolimits}

\def\Ext{\mathop{\rm Ext}\nolimits}
\def\EExt{\mathop{\mathcal Ext}\nolimits}
\def\Tor{\mathop{\rm Tor}\nolimits}

\def\Pic{\mathop{\rm Pic}\nolimits}

\def\Aut{\mathop{\rm Aut}\nolimits}

\def\imm{\mathop{\rm im}\nolimits}

\def\spec{\mathop{\rm spec}\nolimits}

\def\lra{\longrightarrow}

\def\sigg{\mathop{\hbox{$\displaystyle\sum$}}\limits}

\def\hfl#1#2{\smash{\mathop{\ \hbox to 12mm{\rightarrowfill}}
\limits^{\scriptstyle#1}_{\scriptstyle#2} \ }}
\def\hflb#1#2{\smash{\mathop{\hbox to 12mm{\leftarrowfill}}
\limits^{\scriptstyle#1}_{\scriptstyle#2}}}

\def\m#1{{\hbox{$#1$}}}

\def\ot{\otimes}
\def\og{\leavevmode\raise.3ex\hbox{$\scriptscriptstyle\langle\!\langle$}}
\def\fg{\leavevmode\raise.3ex\hbox{$\scriptscriptstyle\,\rangle\!\rangle$}}

\def\nsp{\lbrace 0\rbrace}

\def\dsp{\displaystyle}
\def\Ssect#1#2{\pagebreak[3]\begin{sub}\label{#2}{\sc\small\small
#1}\rm\medskip}

\def\sepsec{\vskip 1.4cm}
\def\sepsub{\vskip 0.7cm}
\def\sepsubsub{\vskip 0.5cm}
\def\sepprop{\vskip 0.5cm}

\def\xmat#1{\[\xymatrix{#1}\]}
\def\flinc{\ar@{^{(}->}}
\def\fleq{\ar@{=}}
\def\flon{\ar@{->>}}
\def\fmaps{\ar@{|-{>}}}
\def\fflat{\ar@{-}}
\def\fimpl{\ar@{=>}}
\def\Nligne{\hfil\break}

\def\ED{\vskip 1cm\end{document}}

\def\BS#1{{\boldsymbol{#1}}}

\def\mm{\mathfrak{m}}

\newcommand{\M}{{\mathbb M}}
\newcommand{\A}{{\mathbb A}}

\newcommand{\Z}{{\mathbb Z}}

\newcommand{\C}{{\mathbb C}}
\newcommand{\Q}{{\mathbb Q}}
\renewcommand{\P}{{\mathbb P}}
\newcommand{\D}{{\mathbb D}}
\newcommand{\F}{{\mathbb F}}
\newcommand{\E}{{\mathbb E}}

\newcommand{\I}{{\mathbb I}}
\newcommand{\J}{{\mathbb J}}

\renewcommand{\L}{{\mathbb L}}

\newcommand{\ka}{{\mathcal A}}

\newcommand{\kd}{{\mathcal D}}
\newcommand{\ke}{{\mathcal E}}
\newcommand{\kf}{{\mathcal F}}
\newcommand{\kg}{{\mathcal G}}
\newcommand{\kh}{{\mathcal H}}
\newcommand{\ki}{{\mathcal I}}
\newcommand{\kj}{{\mathcal J}}
\newcommand{\kk}{{\mathcal K}}

\newcommand{\kn}{{\mathcal N}}
\newcommand{\ko}{{\mathcal O}}

\begin{document}

\def\refname{References}
\def\contentsname{Summary}
\def\proofname{Proof}
\def\abstractname{Resume}

\author{Jean--Marc Dr\'{e}zet}
\address{
Institut de Math\'ematiques de Jussieu - Paris Rive Gauche\\
Case 247\\
4 place Jussieu\\
F-75252 Paris, France}
\email{jean-marc.drezet@wanadoo.fr}
\title[{Coherent sheaves}] {Coherent sheaves on primitive multiple schemes}

\begin{abstract}
A primitive multiple scheme is a Cohen-Macaulay scheme $Y$ such that the 
associated reduced scheme $X=Y_{red}$ is smooth, irreducible, and that $Y$ can 
be locally embedded in a smooth variety of dimension \m{\dim(X)+1}. If $n$ is 
the multiplicity of $Y$, there is a canonical filtration \m{ X=X_1\subset
X_2\subset\cdots\subset X_n=Y}, such that \m{X_i} is a primitive multiple 
scheme of multiplicity $i$. The ideal sheaf \m{\ki_X} of $X$ is a line bundle 
on \m{X_{n-1}} and \m{L=\ki_X/\ki_X^2} is a line bundle on $X$, called the {\em 
associated line bundle} of $Y$.

Even if $X$ is projective, $Y$ needs not to be quasi projective. We define in 
every case the {\em reduced Hilbert polynomial} \m{P_{red,\ko_X(1)}(E)} of a 
coherent sheaf $E$ on $Y$, depending on the choice of an ample line bundle 
\m{\ko_X(1)} on $X$. If $\ke$ is a flat family of sheaves on $Y$ parameterized 
by a smooth curve $C$, then \m{P_{red,\ko_X(1)}(\ke_c)} does not depend on 
\m{c\in C}. We study flat families of sheaves in two important cases: 
the families of {\em quasi locally free sheaves}, and if \m{n=2} those 
of {\em balanced sheaves}. Balanced sheaves are generalizations of vector 
bundles on $Y$, and could be used to expand already known moduli spaces of 
vector bundles on $Y$. 

When $X$ is a smooth projective surface, and $Y$ is of multiplicity 2 we study 
the simplest examples of balanced sheaves: the sheaves $\ke$ such that there is 
an exact sequence
\[0\lra\ki_P\ot L\lra\ke\lra\ki_P=\ke_{|X}\lra 0 \ , \]
where \m{\ki_P\subset\ko_X} is the ideal sheaf of a point \m{P\in X}. They can 
also be described as the ideal sheaves $\ke$ of subschemes of $Y$ supported 
on $P$, and such that \m{\ke_P} is generated by two elements whose images in 
\m{\ko_{X,P}} generate the maximal ideal. There is a moduli space for such 
sheaves, which is an affine bundle on $X$ with associated vector bundle 
\m{T_X\ot L} (where \m{T_X} is the tangent bundle of $X$). The associated class 
in \m{H^1(X,T_X\ot L)} can be determined.
\end{abstract}

\maketitle
\tableofcontents


\vskip 1cm

\section{Introduction}\label{intro}

A {\em primitive multiple scheme} is a Cohen-Macaulay scheme $Y$ over $\C$ such 
that:
\begin{enumerate}
\item[--] $Y_{red}=X$ is a smooth connected variety, 
\item[--] for every closed point \m{x\in X}, there exists a neighborhood $U$ 
of $x$ in $Y$, and a smooth variety $S$ of dimension \ \m{\dim(X)+1} such that 
$U$ is isomorphic to a closed subscheme of $S$.
\end{enumerate}
We call $X$ the {\em support} of $Y$. It may happen that $Y$ is 
quasi-projective, and in this case it is projective if $X$ is.

For every closed subscheme \m{Z\subset Y}, let \m{\ki_Z} (or \m{\ki_{Z,Y})} 
denote the ideal sheaf of $Z$ in $Y$. For every positive integer $i$, 
let \m{X_i} be the closed subscheme of $Y$ corresponding to the ideal sheaf 
\m{\ki_X^i}. The smallest integer $n$ such that \ \m{X_n=Y} \ is called the 
{\em multiplicity} of $Y$. For \m{1\leq i\leq n}, \m{X_i} is a primitive 
multiple scheme of multiplicity $i$, \m{L=\ki_X/\ki_{X_2}} is a line bundle on 
$X$, and we have \ \m{\ki_{X_{i}}/\ki_{X_{i+1}}=L^i}. 
We call $L$ the line bundle on $X$ {\em associated} to $Y$. The ideal sheaf 
\m{\ki_{X}} can be viewed as a line bundle on \m{X_{n-1}}. If \m{n=2}, $Y$ is 
called a {\em primitive double scheme}.

The simplest case is when $Y$ is contained in a smooth variety $S$ of dimension 
\m{\dim(X)+1}. Suppose that $Y$ has multiplicity $n$. Let \m{P\in X} and 
\m{f\in\ko_{S,P}}  a local equation of $X$. Then we have \ 
\m{\ki_{X_i,P}=(f^{i})} \ for \m{1<j\leq n} in $S$, in particular 
\m{\ki_{Y,P}=(f^n)} (where \m{\ki_{X_i}, \ki_Y} are the ideal sheaves of 
\m{X_i} and \m{Y} in $S$ respectively), and \ \m{L=\ko_X(-X)} .

For any \m{L\in\Pic(X)}, the {\em trivial primitive scheme} of multiplicity 
$n$, with induced smooth variety $X$ and associated line bundle $L$ on $X$ is 
the $n$-th infinitesimal neighborhood of $X$, embedded by the zero section in 
the dual bundle $L^*$, seen as a smooth variety.

The primitive multiple curves where defined in \cite{fe}, \cite{ba_fo}. 
Primitive double curves were parameterized and studied in \cite{ba_ei} and 
\cite{ei_gr}. More results on primitive multiple curves can be found in 
\cite{dr2}, \cite{dr1}, \cite{dr4}, \cite{dr5}, \cite{dr6}, \cite{dr7}, 
\cite{dr8}, \cite{dr9}, \cite{ch-ka}, \cite{sa1}, \cite{sa-vi}, \cite{sa3}. 
Some 
primitive double schemes are studied in \cite{b_m_r}, \cite{ga-go-pu} and 
\cite{gonz1}. The case of varieties of any dimension is studied in \cite{dr10}, 
where the following subjects were treated:
\begin{enumerate}
\item[--] construction and parameterization of primitive multiple schemes,
\item[--] obstructions to the extension of a vector bundle on \m{X_m} to 
\m{X_{m+1}},
\item[--] obstructions to the extension of a primitive multiple scheme of 
multiplicity $n$ to one of multiplicity \m{n+1}.
\end{enumerate}
In \cite{dr11}, the construction and properties of fine moduli spaces of vector 
bundles on primitive multiple schemes are described. Suppose that \m{Y=X_n} is 
of multiplicity $n$, and can be extended to \m{X_{n+1}} of multiplicity 
\m{n+1}, and let \m{M_n} be a fine moduli space of vector bundles on 
\m{X_n}. With suitable hypotheses, a fine moduli space \m{M_{n+1}} for the 
vector bundles on \m{X_{n+1}} whose restriction to \m{X_n} belongs to 
\m{M_n} can be constructed. This applies in particular to Picard groups.

The subject of this paper is the study of more general coherent sheaves on 
primitive multiple schemes. More precisely:
\begin{enumerate}
\item[--] We will define an invariant of coherent sheaves on $Y$ depending on 
the choice of an ample line bundle on $X$, the {\em reduced Hilbert 
polynomial}, which is the same as the usual Hilbert polynomial when $Y$ is 
projective. But it is defined even when $Y$ is not quasi-projective.
\item[--] We define and study {\em primitive multiple rings}. Among them are 
the local rings of closed points in a primitive multiple scheme. They make the 
presentations more readable.
\item[--] We give some properties of {\em quasi locally free sheaves} already 
studied on primitive multiple curves, and of {\em balanced sheaves}, which are 
natural generalizations of vector bundles. The main subjects are flat families 
of such sheaves.
\item[--] We study the simplest balanced sheaves on \m{X_2} which are not 
locally free, i.e. sheaves whose restriction to $X$ is the ideal sheaf of a 
point.
\end{enumerate}

\sepsub

\Ssect{Reduced Hilbert polynomials}{RHP}

It may happen that $Y$ is not quasi-projective. For example, there are exactly 
two non trivial primitive multiple schemes such that $X$ is a projective space 
of dimension \m{\geq 2}: the two with \m{X=\P_2}, one of multiplicity 2, the 
other of multiplicity 4 (cf. \cite{dr10}). On these two schemes the only line 
bundle is the trivial one, so they are not quasi-projective, and there is no 
notion of Hilbert polynomial for sheaves on such schemes.

If $Z$ is a projective scheme,  \m{\ko_Z(1)} an ample line bundle  and $E$ a 
coherent sheaf on $Z$, let \m{P_{\ko_Z(1)}(E)} be the Hilbert polynomial of $E$ 
with respect to \m{\ko_Z(1)}.

Let \ \m{Y=X_n} \ be a primitive multiple scheme of multiplicity $n$ and 
associated smooth scheme $X$. Suppose that $X$ is projective. Let \m{\ko_X(1)} 
be an ample line bundle on $X$. Let $\ke$ be a coherent sheaf on $Y$. There 
exist filtrations
\begin{equation}\label{equ23}
\kf_m=0\subset\kf_{m-1}\subset\cdots\subset\kf_1\subset\kf_0=\ke
\end{equation}
such that, for \m{0\leq i<m}, \m{\kf_i/\kf_{i+1}} is 
supported on $X$ (for example the two canonical filtrations, cf. \ref{QLL}). We 
will see (proposition \ref{prop1}) that
\[P_{red,\ko_X(1)}(\ke) = \ \sigg_{i=0}^{m-1}P_{\ko_X(1)}(\kf_i/\kf_{i+1})\]
does not depend on the filtration $(\ref{equ23})$. It is called the {\em 
reduced Hilbert polynomial} of $\ke$ with respect to \m{\ko_X(1)}. This expands 
to higher dimensions the notions of generalized rank and degree defined in 
\cite{dr2}, \cite{dr4} on primitive multiple curves. A similar definition can 
be made for {\em reduced Chern classes}. These invariants behave correctly with 
exact sequences of sheaves.

If \m{\ko_X(1)} can be extended to a line bundle \m{\ko_Y(1)} on $Y$, then 
\m{\ko_Y(1)} is ample and \Nligne \m{P_{red,\ko_X(1)}=P_{\ko_Y(1)}}.

The reduced Hilbert polynomial of a coherent sheaf is invariant by flat 
deformations of the sheaf: let $C$ be an irreducible smooth curve, and $\ke$ a 
coherent sheaf on \m{X_n\times C}, flat on $C$. Then the map 
\[\xymatrix@R=5pt{C\ar[r] & \Q[T]\\ c\fmaps[r] & P_{red,\ko_X(1)}(\ke_c)
}\]
is constant (corollary \ref{conj}). To prove this we need a preliminary result 
that was pointed out by J\'anos Koll\'ar, which has its own interest: suppose 
that we have a filtration
\[0=\kg_m\subset\kg_{m-1}\subset\cdots\subset\kg_1\subset\kg_0=\ke\]
such that for \ \m{0\leq i<m}, \m{\kg_i/\kg_{i+1}} is supported on 
\m{X\times C}. Then there exists a filtration
\[0=\kf_m\subset\kf_{m-1}\subset\cdots\subset\kf_1\subset\kf_0=\ke\] such that
the \m{\kf_i} are flat on $C$, as well as the quotients \m{\kf_i/\kf_{i+1}}, 
which are supported on \m{X\times C}, and that there is a finite subset 
\m{\Sigma\subset C} such that the two filtrations coincide on \ 
\m{X_n\times(C\backslash\Sigma)} (theorem \ref{theo8}).
\end{sub}

\sepsub

\Ssect{Quasi locally free sheaves}{QLL}

\begin{subsub}\label{canf}Canonical filtrations -- \rm Let $\ke$ be a coherent 
sheaf on $Y$. The {\em first canonical filtration}  of $\ke$ is 
\[\ke_n=0\subset \ke_{n-1}\subset\cdots\subset \ke_{1}\subset \ke_0=\ke \ , \]
where for \m{0\leq i< n}, \m{\ke_{i+1}} is the kernel of the canonical
surjective morphism \ \m{\ke_i\to\ke_{i\mid X}}. For \m{0\leq i<n}, let \
\m{G_i(\ke)=\ke_{i}/\ke_{i+1}=\ke_{i\mid X}} . The {\em second canonical 
filtration of $\ke$} is
\[\ke^{(0)}=\nsp\subset \ke^{(1)}\subset\cdots\subset
\ke^{(n-1)}\subset \ke^{(n)}=\ke \ , \]
where for \m{1\leq i\leq n-1}, \m{\ke^{(i)}} is the maximal subsheaf 
annihilated by \m{\ki_X^i}. For \m{0<i\leq n}, let \ \m{G^{(i)}(\ke)= 
\ke^{(i)}/\ke^{(i-1)}} . The sheaves \m{G^{(i)}(\kf)} and \m{G_i(\kf)} are 
supported on $X$.
\end{subsub}

\sepprop

Let $\ke$ be a coherent sheaf on $Y$ and \m{x\in X} a closed point. We say that 
$\ke$ is {\em quasi locally free at $x$} if there exist integers \m{m_i\geq 0}, 
\m{1\leq i\leq n}, such that there is an isomorphism of \m{\ko_{Y,x}}-modules
\[\ke_x \ \simeq \ \bigoplus_{1\leq i\leq n}\ko_{X_i,x}^{\oplus m_i} \ . \] 
The set $U$ of such points $x$ is nonempty and open, the integers \m{m_i} are 
uniquely determined and do not depend on $x$. If \m{U=X}, $\ke$ is called {\em 
quasi locally free}. The sequence \m{(m_1,\ldots,m_n)} is called the {\em type} 
of $\ke$. If $\kf$ is a coherent sheaf on $Y$, $\kf$ is torsion free (resp. 
quasi locally free) if and only if all the \m{G^{(i)}(\kf)} (resp. \m{G_i(\kf)})
are torsion free (resp. locally free).

We will see (in theorem \ref{theo4}) that if we have a flat family of quasi 
locally free sheaves $\ke$ on $Y$, all of the same type, then the \m{\ke_i} and 
\m{\ke^{(i)}} form flat families of sheaves. Flat families of such sheaves 
having this property were called {\em good families} in \cite{dr10b}, 3.2. 

\end{sub}

\sepsub

\Ssect{Balanced sheaves}{BAS}

Let $\E$ be a vector bundle on \m{Y=X_n} and \ \m{E=\E_{|X}}. Then for \m{0\leq 
i\leq n}, \m{\E_i} is locally free on \m{X_{n-i}}, \m{\E_i=\E^{(n-i)}} and \ 
\m{G_i(\E)=G^{(n-i)}(\E)=E\ot L^i}.

Let $\ke$ be a coherent sheaf on \m{X_n}. We say that $\ke$ is {\em balanced} 
if \ \m{\ke_i=\ke^{(n-i)}} \ for \m{1\leq i\leq n} (see \ref{ext_hi} for an 
equivalent definition). Here also \ \m{G_i(\ke)=G^{(n-i)}(\ke)=\ke_{|X}\ot L^i}.

Balanced sheaves are the most natural sheaves that could be extensions of 
torsion free sheaves on $X$ to \m{X_n}. To build moduli spaces of these 
sheaves, the work of \cite{dr11} on vector bundles should be extended to 
torsion free sheaves.

\sepprop

\begin{subsub} The case of primitive double schemes -- \rm
Suppose that \m{n=2}, i.e. $Y$ is a primitive double scheme. We prove that the 
fact that a sheaf is balanced is an {\em open property}: let $C$ be a smooth 
curve and $\ke$ a family of coherent sheaves on $Y$, parameterized by $C$, flat 
on $C$. Suppose that for some closed point \m{c_0\in C}, \m{\ke_{c_0}} is 
balanced. Then there exists an open neighborhood $U$ of \m{c_0} such that for 
every \m{c\in U}, \m{\ke_c} is balanced (cf. theorem \ref{theo5}).

Let $\ke$ a flat family of balanced sheaves on \m{X_2} parameterized by a 
smooth connected curve $C$. We prove that the restrictions \m{\ke_{c|X}} form a 
flat family of sheaves on $X$ parameterized cy $C$ (corollary \ref{cor6}).
\end{subsub}

\sepprop

\begin{subsub} Extension of sheaves on a smooth surface to multiplicity 2 -- 
\rm 
Suppose that \m{n=2} and \ \m{\dim(X)=2}. Let $S$ be a smooth variety and 
$\boldsymbol{F}$ a coherent sheaf on \m{X\times S} such that
\begin{enumerate}
\item[--] $\boldsymbol{F}$ is flat on $S$.
\item[--] For every \m{s\in S}, \m{\boldsymbol{F}_s} is torsion free and simple.
\item[--] $\dim(\Ext^1_{\ko_X}({\bf F}_s,{\bf F}_s\ot L))$ is independent of 
$s\in S$.
\end{enumerate}
We prove in \ref{fam_bal} that the set of points \m{s\in S} such that 
\m{\boldsymbol{F}_s} can be extended to a balanced sheaf on \m{X_2} is closed 
in $S$. This could allow to extend some results on moduli spaces of vector 
bundles on \m{X_2} (cf. \cite{dr11}) to moduli spaces of balanced sheaves.
\end{subsub}

\end{sub}

\sepsub

\Ssect{Ideal sheaves}{IS}

Let $\boldsymbol{\ki}$ be the ideal sheaf of a zero dimensional subscheme $Z$
of \m{X_n}. Suppose that for every closed point $x$ of $Z$, the 
\m{\ko_{X_n,x}}-module \m{\boldsymbol{\ki}_x} is generated by a regular 
sequence. This is equivalent to the fact that the image of 
\m{\boldsymbol{\ki}_x} in \m{\ko_{X,x}} is generated by a regular sequence. 
Then $\boldsymbol{\ki}$ is balanced (cf. propositions \ref{prop20} and 
\ref{prop21}).

\sepprop

\begin{subsub} Ideal sheaves on primitive double surfaces -- \rm
We suppose now that \ \m{\dim(X)=2} \ and \m{n=2}. Let \m{P\in X} and \m{\ki} 
an ideal sheaf on \m{X_2} of a subscheme whose support is \m{\{P\}}. Suppose 
that \m{\ki_P} is generated by two elements \m{x,y} whose images in 
\m{\ko_{X,P}} are generators of the maximal ideal. We show that
\begin{enumerate}
\item[--] The set of such ideal sheaves, with fixed $P$, has a natural 
structure of affine space with associated vector space \ \m{T_{X,P}\ot L_P}, 
where \m{T_X} is the tangent bundle of $X$ (cf. proposition \ref{prop24}).
\item[--] The set $\I$ of all these sheaves, for all $P$, has a natural 
structure of an affine bundle on \m{X}, with associated vector bundle \m{T_X\ot 
L}.
\item[--] We build in \ref{uni_s} a {\em universal family of ideal sheaves} 
$\boldsymbol{\kj}$ parameterized by $\I$, such that for every ideal sheaf 
\m{I\in\I}, \m{\boldsymbol{\kj}_I\simeq I}.
\end{enumerate}
We have a canonical exact sequence of vector bundles on $X$:
\[0\lra L=\ki_X\lra\Omega_{X_2|X}\lra\Omega_X\lra 0 \ , \]
associated with \ \m{\sigma\in\Ext^1_{\ko_X}(\Omega_X,L)=H^1(X,T_X\ot L)}. To 
the affine bundle $\I$ is also associated \ \m{\eta\in H^1(X,T_X\ot L)}. We 
prove that \ \m{\C\sigma=\C\eta}.

To obtain a moduli space of sheaves we consider the set \ 
\m{\boldsymbol{I}(X_2)=\I\times\Pic(X_2)} \ of sheaves of the form 
\m{I\ot D}, where \m{I\in\I} and \m{D\in\Pic(X_2)}. Let 
\m{p_1:\boldsymbol{I}(X_2)\to\I}, \m{p_1:\boldsymbol{I}(X_2)\to\Pic(X_2)} \ be 
the projections. The exists an open cover \m{(P_i)_{i\in I}} of \m{\Pic(X_2)} 
such that for every \m{i\in I}, there is a Poincaré bundle \m{\kd_i} on 
\m{X_2\times P_i}. This defines \m{\Pic(X_2)} as a {\em fine moduli space} of 
vector bundles (cf. \cite{dr11}, \ref{fin-mod}, \ref{pic-sch}).

Let \ \m{\boldsymbol{\kh_i}=p_1^\sharp(\boldsymbol{\kj})\ot 
p_2^\sharp(\kd_i)}, which is a sheaf on \ \m{\I\times P_i}. For  every ideal 
$\ki$ in $\I$ and \ \m{D\in P_i}, \m{\boldsymbol{\kh}_{\ki,D}=\ki\ot\kd_{i,D}
\simeq\ki\ot D} . Then \m{\boldsymbol{I}(X_2)}, with the open cover 
\m{(\I\times P_i)_{i\in I}} and the sheaves \m{\boldsymbol{\kh_i}} is a fine 
moduli space in the sense of \ref{fin-mod} (proposition \ref{prop35}).
\end{subsub}

\end{sub}

\sepsub

\Ssect{Outline of the paper}{intro-5}

In Chapter 2, we give several preliminary definitions and technical results 
that will be used in the next chapters.

In Chapter 3, we define and study the primitive multiple rings. Some results 
obtained here are used in the next chapters.

In Chapter 4 we recall the definitions and properties of primitive multiple 
schemes, and we prove the results obtained on quasi locally free sheaves.

In Chapter 5 we define and give some properties of reduced Hilbert polynomials
of coherent sheaves on primitive multiple schemes. 

The Chapter 6 is devoted to balanced sheaves on primitive multiple schemes.

In Chapter 7 we study the balanced sheaves on a primitive double surface 
\m{X_2}, whose restriction to $X$ are ideal sheaves of a point.

\end{sub}

\sepprop

{\bf Notations and terminology: } -- Let \m{x_0,\ldots,x_k\in\C}, not all zero. 
We will also denote by \m{(x_0,\ldots,x_k)} the element \m{\C.(x_0,\ldots,x_k)} 
of \m{\P_k}.

-- A {\em scheme} is a noetherian separated scheme over $\C$.

-- If $X$ is a scheme and \m{Y\subset X} a subscheme, \m{\ki_{Y,X}} (or 
\m{\ki_Y}) denotes the ideal sheaf of $Y$ in $X$. 

-- Let $X$, $Y$, $Z$ be schemes, $\ke$ a coherent sheaf on \m{X\times Z}, and \ 
\m{f:Y\to Z} \ a morphism. Then \ \m{f^\sharp(\ke)=(I_X\times f)^*(\ke)}.

\sepsec

\section{Preliminaries}\label{prelim}

\Ssect{Canonical class of a line bundle}{can_class}

Let $Z$ be a scheme over $\C$ and $L$ a line bundle on $Z$. To $L$ one 
associates an element \m{\nabla_0(L)} of \m{H^1(Z,\Omega_Z)}, called the {\em 
canonical class of $L$}. If $Z$ is smooth and projective, and \m{L=\ko_Z(Y)}, 
where \m{Y\subset Z} is a smooth hypersurface, then \m{\nabla_0(L)} is the 
cohomology class of $Y$.

Let \m{(Z_i)_{i\in I}} be an open cover of $Z$ such that $L$ is defined by a 
cocycle \m{(\theta_{ij})}, \m{\theta_{ij}\in\ko_Z(Z_{ij})^*}. Then 
\m{\dsp\Big(\frac{d\theta_{ij}}{\theta_{ij}}\Big)} is a cocycle which 
represents \m{\nabla_0(L)}.
\end{sub}

\Ssect{Extensions of modules}{extmod}

Let $R$ be a commutative ring and $M$ a $R$-module. Suppose that we have a free 
resolution
\xmat{\cdots\F_2\ar[r]^-{\phi_2} & \F_1\ar[r]^-{\phi_1} & \F_0\ar[r]\flon[r] & 
M \ . }
Let $N$ be a $R$-module. We have exact sequences
\[\Hom(\F_0,N)\lra\Hom(\imm(\phi_1),N)\hfl{\gamma}{}\Ext^1_R(M,N)\lra 0 \ , \]
\[0\lra\Hom(\imm(\phi_1),N)\hfl{\beta}{}\Hom(\F_1,N)\lra\Hom(\F_2,N) \ . \]
Let \ \m{0\to N\to P\to M\to 0} \ be an exact sequence, corresponding to \ 
\m{\sigma\in\Ext^1_R(M,N)}. Suppose that \m{\sigma=\gamma(f)}, 
\m{f:\imm(\phi_1)\to N}, and \m{f_1=\beta(f)}. Then $P$ is constructed as 
follows: let \ \m{\eta=f_1\oplus\phi_1:\F_1\to N\oplus\F_0}. Then \ 
\m{P=\coker(\eta)}. The inclusion \ \m{N\hookrightarrow P} is induced by \ 
\m{N\hookrightarrow N\oplus\F_0}, and \ \m{P\twoheadrightarrow M} \ by \ 
\m{\F_0\twoheadrightarrow M}.

\end{sub}

\sepsub

\Ssect{Refinements of filtrations}{RF}

Let $Z$ be a scheme, $\ke$ a sheaf of \m{\ko_Z}-modules, and two filtrations of
$\ke$
\[\kd : \ \kd_n=0\subset\kd_{n-1}\subset\cdots\subset\kd_1\subset\kd_0=\ke
\ , \]
\[\kf : \ \kf_m=0\subset\kf_{m-1}\subset\cdots\subset\kf_1\subset\kf_0=\ke
\ . \]
We say that $\kd$ is a {\em refinement} of $\kf$ if for every $j$, \m{0\leq 
j\leq m}, there exists $i$, \m{0\leq i\leq n}, such that \ \m{\kd_i=\kf_j}. In 
this case we write \ \m{\kd\to\kf}.

Two filtrations \m{(\kd_i)_{0\leq i\leq n}}, \m{(\kf_j)_{0\leq j\leq m}} of
$\ke$ are called {\em similar} if \m{m=n}, and if there exists a permutation
$\sigma$ of \m{\{0,\ldots,n\}} such that for every \m{i\in\{0,\ldots,n\}} we
have\ \m{\kd_i/\kd_{i+1}\simeq\kf_{\sigma(i)}/\kf_{\sigma(i)+1}}. In this case 
we write \ \m{\kd\simeq\kf}.

\sepprop

\begin{subsub}\label{theo7}{\bf Proposition: } Let
\[\kd \ : \ \kd_n=0\subset\kd_{n-1}\subset\cdots\subset\kd_1\subset\kd_0=\ke \ 
, \]
\[\kf \ : \ \kf_m=0\subset\kf_{m-1}\subset\cdots\subset\kf_1\subset\kf_0=\ke \ 
 \]
be two filtrations of $\ke$. Then there exist refinements \m{\kd'} of $\kd$, 
\m{\kf'} of $\kf$, which are similar.
\end{subsub}

\begin{proof} {\bf Step 1 -- } We first suppose that \m{n=2}. For 
\m{(\kd'_k)_{0\leq k\leq p}} we take
\[0=\kd_1\cap\kf_m\subset\kd_1\cap\kf_{m-1}\subset
\cdots\kd_1\cap\kf_1\subset\kd_1\subset\kd_1+\kf_{m-1}
\subset\cdots\subset\kd_1+\kf_1\subset\kd_1+\kf_0=\ke \ ,
\]
and we build \m{(\kf'_k)_{0\leq k\leq p}} by inserting \
\m{\kf_{j+1}+\kd_1\cap\kf_j} \ between \m{\kf_{j+1}} and \m{\kf_j}.
We have
\[(\kf_{j+1}+\kd_1\cap\kf_j)/\kf_{j+1} \ \simeq \
(\kd_1\cap\kf_j)/(\kd_1\cap\kf_{j+1}) \ , \]
\[\kf_j/(\kf_{j+1}+\kd_1\cap\kf_j) \ \simeq \
(\kd_1+\kf_j)/(\kd_1+\kf_{j+1}) \ . \]

{\bf Step 2 -- } Now suppose that we have two similar filtrations of $\ke$
\[\kg \ : \ \kg_p=0\subset\kg_{p-1}\subset\cdots\subset\kg_1\subset\kg_0=\ke \ 
, \]
\[\kh : \ \kh_p=0\subset\kh_{p-1}\subset\cdots\subset\kh_1\subset\kh_0=\ke \ , 
\]
and that we have a refinement \m{\kg'} of $\kg$. Then there exists a refinement 
\m{\kh'} of $\kh$ which is similar to \m{\kg'}: let $\sigma$ be a permutation 
of \m{\{0,\ldots,p\}} such that \ \m{\kg_i/\kg_{i+1}\simeq\kh_{\sigma(i)}/
\kh_{\sigma(i)+1}}. Then to each \m{\kg'_k} such that \ 
\m{\kg_{i+1}\subset\kg'_k\subset\kg_i} corresponds a sheaf \m{\kk} such that \
\m{\kh_{\sigma(i)+1}\subset\kk\subset\kh_{\sigma(i)}}. All these sheaves make 
the filtration \m{\kh'}.

{\bf Step 3 -- } Now we prove proposition \ref{theo7} by induction on $n$. It 
is already proved for \m{n=2}. Suppose that it is true for \m{n-1\geq 2}. Then 
we consider the filtrations
\[\kd^0 \ : \ 
\kd^0_{n-1}=0\subset\kd^0_{n-2}=\kd_{n-2}\subset\cdots\subset\kd^0_1=\kd_1
\subset\kd^0_0=\ke\]
and $\kf$. There are refinements \m{\kd''=(\kd''_k)_{0\leq k\leq p}} of 
\m{\kd^0} and \m{\kf''} of $\kf$ which are similar. There exists an integer $q$ 
such that \m{\kd''_q=\kd_{n-2}}. We have filtrations of \m{\kd_{n-2}}
\[0\subset\kd_{n-1}\subset\kd_{n-2} \ , \]
\[0=\kd''_p\subset\kd''_{p-1}\subset\cdots\subset\kd''_q=\kd_{n-2} \ . \]
From step 1, we can find refinements \m{(\kd^1_k)_{0\leq k\leq r}} of the 
first, \m{(\kd^2_k)_{0\leq k\leq r}} of the second, which are similar. Now we 
build two refinements \m{\kd^3}, \m{\kd^4}, of \m{\kd^0}, the first by 
replacing in \m{\kd''}, \m{(\kd''_p,\ldots,\kd''_q)} with 
\m{(\kd^1_r,\ldots,\kd^1_0)}, and the second by replacing in \m{\kd''}, 
\m{(\kd''_p,\ldots,\kd''_q)} with \m{(\kd^2_r,\ldots,\kd^2_0)}. We have \ 
\m{\kd^3\simeq\kd^4}, \m{\kd^3\to\kd}, \m{\kd^4\to\kd''} and 
\m{\kd''\simeq\kf''}.

From step 2 we can find a refinement \m{\kf^4} of \m{\kf''} such that 
\m{\kf^4\simeq\kd^4}. So we can take \m{\kd'=\kd^3}, \m{\kf'=\kf^4}.
\end{proof}

\end{sub}

\sepsub

\Ssect{Affine bundles}{aff_b}

Let \ \m{f:\ka\to S} \ be a morphism of schemes, and \m{r\geq 0} an integer. We 
say that $f$ (or $\ka$) is an {\em affine bundle} of rank $r$ over $S$ if there 
exists an open cover \m{(S_i)_{i\in I}} of $S$ such that for every \m{i\in I} 
there is an isomorphism \ \m{\tau_i:f^{-1}(S_i)\to S_i\times\C^r} \ over 
\m{S_i} such that for every distinct \m{i,j\in I}, \m{\tau_j\circ 
\tau_i^{-1}:S_{ij}\times\C^r\to S_{ij}\times\C^r} \ is of the form
\[(x,u)\longmapsto (x,A_{ij}(x)u+b_{ij}(x)) \ , \]
where \m{A_{ij}} is an \m{r\times r}-matrix of elements of \m{\ko_S(S_{ij})} 
and \m{b_{ij}} is a morphism from \m{S_{ij}} to \m{\C^r}. We have then the 
cocycle relations
\[A_{ij}A_{jk} \ = \ A_{ik} \ , \quad b_{ik} \ = \ A_{ij}b_{jk}+b_{ij} \ . \]
The first relation shows that the family \m{(A_{ij})} defines a vector bundle 
$\A$ on $S$, and the second that \m{(b_{ij})} defines \ \m{\lambda\in 
H^1(S,\A)} (according to \cite{dr10}, 2.1.1). 

The vector bundle $\A$ is uniquely defined, as well as \ 
\m{\eta(\ka)=\C\lambda\in\big(\P(H^1(S,\A))\cup\{0\}\big)/\Aut(\A)}.

We say that {\em $\ka$ is a vector bundle} if $f$ has a section. This is the 
case if and only \m{\lambda=0}, and then a section of $f$ induces an 
isomorphism \m{\ka\simeq\A} over $S$.

For every closed point \m{s\in S} there is a canonical action of the additive 
group \m{\A_s} on \m{\ka_s}
\[\xymatrix@R=5pt{\A_s\times\ka_s\ar[r] & \ka_s\\ (u,a)\fmaps[r] & a+u}\]
such that for every \m{a\in\ka_s}, $u\mapsto a+u$ is an isomorphism 
\m{\A_s\simeq\ka_s}.

\end{sub}

\sepsub

\Ssect{Locally free resolutions}{loc_free}

Let $X$, $S$ be smooth algebraic varieties, with $X$ projective. Let $\ke$ be a 
flat family of torsion free sheaves on $X$ parameterized by $S$, and
\[\cdots\E_m\hfl{f_m}{}\E_{m-1}\hfl{f_{m-1}}{}\cdots\lra\E_0\hfl{f_0}{}\ke\lra 
0\]
an exact sequence, where the \m{\E_i} are flat families of torsion free sheaves 
on $X$ parameterized by $S$.

\sepprop

\begin{subsub}{\bf Lemma: }\label{lem7} The sheaf $\ke$ is torsion free.
\end{subsub}
\begin{proof} Let \m{(x,s)\in X\times S}, \m{e\in\ke_{(x,s)}}, 
\m{\alpha\in\ko_{X\times S,(x,s)}} be such that \m{\alpha e=0}. Let $U$ be an 
open neighborhood of \m{(x,s)} such that \m{e\in\ke(U)}, 
\m{\alpha\in\ko_{X\times S}(U)}. Suppose that \m{\alpha\not=0}. If 
\m{(x',s')\in U}, and \m{\ke_{s',x'}} is a free 
\m{\ko_{X\times\{s'\},x'}}-module , then \m{\ke_{(x',s')}} is a free 
\m{\ko_{X\times S,(x',s')}}-module. Hence \m{e=0} on a neighborhood of 
\m{(x',s')}. If \m{s'\in S} is such that \ 
\m{U_{s'}=U\cap(X\times\{s'\})\not=\emptyset}, then \m{\ke_{|U_{s'}}} is 
locally free on a nonempty open subset of \m{U_{s'}}. Hence \m{e_{|U_{s'}}=0}. 
Since this is true for every such \m{s'}, we have \m{e=0}.
\end{proof}

\sepprop

\begin{subsub}{\bf Lemma: }\label{lem6} Let \m{H\subset S} be a hypersurface, 
and $\kf$ a torsion free sheaf on \m{X\times S}. Then we have
\[\Tor^1_{\ko_{X\times S}}(\kf,\ko_{X\times H}) \ = \ 0 \ .\]
\end{subsub}
\begin{proof}
We have a locally free resolution
\[0\lra\ko_{X\times S}(-(X\times H))\hfl{\gamma}{}\ko_{X\times 
S}\lra\ko_{X\times H}\lra 0\]
that we use to compute \m{\Tor^1_{\ko_{X\times S}}(\ke,\ko_{X\times H})}. It is 
the kernel of
\[I_\kf\ot\gamma:\kf(-(X\times H))\lra\kf \ , \]
which is injective since $\kf$ is torsion free.
\end{proof}

\sepprop

\begin{subsub}{\bf Proposition: }\label{prop12} For every closed point \m{s\in 
S}, the sequence
\[\cdots\E_{m,s}\hfl{f_{m,s}}{}\E_{m-1,s}\hfl{f_{m-1,s}}{}\cdots\lra\E_{0,s}
\hfl{f_{0,s}}{}\ke_s\lra 0\]
is exact.
\end{subsub}
\begin{proof}
For \m{0\leq i<m}, let \ \m{N_i=\ker(f_i)}, so that we have exact sequences
\[0\lra N_i\lra\E_i\lra N_{i-1}\lra 0 \qquad \text{if} \quad i>0 \ ,\]
\[0\lra N_0\lra\E_0\lra\ke\lra 0 \ . \]
Let \m{d=\dim(S)}, and \m{S_1,\ldots,S_d} hypersurfaces containing $s$, smooth 
at $s$ and that intersect transversally at $s$. By replacing $S$ with an open 
neighborhood of $s$ we can assume the \ \m{S_1\cap\cdots\cap S_d=\{s\}}.

From lemma \ref{lem6}, since the \m{N_i} are torsion free, we have exact 
sequences
\[0\lra N_{i|X\times S_1}\lra\E_{i|X\times S_1}\lra N_{i-1|X\times S_1}\lra 0 
\qquad \text{if} \quad i>0 \ ,\]
\[0\lra N_{0|X\times S_1}\lra\E_{0|X\times S_1}\lra\ke_{|X\times S_1}\lra 0 \ , 
\]
i.e. we have a locally free resolution of \m{\ke_{|X\times S_1}}
\[\cdots\E_{m|X\times S_1}\hfl{f_m}{}\E_{m-1|X\times S_1}\hfl{f_{m-1}}{}
\cdots\lra\E_{0|X\times S_1}\hfl{f_0}{}\ke_{|X\times S_1}\lra 0 \ , \]
and by lemma \ref{lem7}, \m{\ke_{|X\times S_1}} and the \m{N_{i|X\times S_1}} 
are torsion free. So we can continue this process, and by induction on $d$ we 
finally arrive at the resolution of proposition \ref{prop12}.
\end{proof}

\end{sub}

\sepsub

\Ssect{Fine moduli spaces of sheaves}{fin-mod}

Let $\chi$ be a set of isomorphism classes of coherent sheaves on a scheme $Z$. 
Suppose that $\chi$ is {\em open}, i.e. for every scheme $V$ and coherent sheaf 
$\E$ on \m{Z\times V}, flat on $V$, if \m{v\in V} is a closed point such that 
\m{\E_v\in\chi}, then there exists an open neighborhood $U$ of $v$ such that 
\m{\E_u\in\chi} for every closed point \m{u\in U}.
A {\em fine moduli space for $\chi$ }is the data of a scheme $M$ and of
\begin{enumerate}
\item[--] a bijection
\[\xymatrix@R=5pt{M^0\ar[r] & \chi\\ m\fmaps[r] & E_m
}\]
(where \m{M^0} denotes the set of closed points of $M$),
\item[--] an open cover \m{(M_i)_{i\in I}} of $M$, and for every \m{i\in I}, a 
coherent sheaf \m{\ke_i} on \m{Z\times M_i} such that for every \m{m\in M_i}, 
\m{\ke_{i,m}\simeq E_m},
\end{enumerate}
such that:
for any scheme $S$, any coherent sheaf $\kf$ on \m{Z\times S}, flat on $S$, 
such that for any closed point \m{s\in S}, \m{\kf_s\in\chi}, there is a 
morphism \ \m{f_\kf:S\to M} \ such that for every \m{s\in S}, if \m{m=f_\kf(s)} 
then \m{\kf_s\simeq E_m}, and if \m{m\in M_i}, then there exists an open 
neighborhood \m{U} of $s$ such that \ \m{f_\kf(U)\subset M_i} and
\ \m{(I_Z\times f_{\kf|U})^*(\ke_i)\simeq\kf_{|Z\times U}}.

\end{sub}

\sepsec

\section{Primitive multiple rings}\label{PMR}

Let $R$ be an entire noetherian ring and $K$ its fraction field.

Let $n$ be a positive integer and \ \m{R[n]=R[t]/(t^n)}, which we call the {\em 
primitive multiple ring of multiplicity} $n$ associated to $R$.

We will also consider rings \m{A,A_n} isomorphic to \m{R,R[n]} respectively, 
with the ideal \m{\ki\subset A_n} corresponding to \m{(t)}, and \m{A=A_n/\ki}. 
Sometimes it is useful not to specify a generator of $\ki$, in particular when 
considering families of rings.

An element $u$ of $R[n]$ can be written in an unique way as
\[u \ = \ \sigg_{i=0}^{n-1}u_it^i \ , \]
with \m{u_i\in R} for \ \m{0\leq i\leq n}. Then $u$ is a zero divisor if and 
only if \m{u_0=0}, i.e. if $u$ is a multiple of $t$. Let \m{S_n} be the 
multiplicative system of non zero divisors of \m{R[n]}, i.e. \ \m{S_n=\{u\in 
R[n];u_0\not=0\}}. It is easy to see that \ \m{S_n^{-1}R[n]=K[t]/(t^n)=K[n]} 
(the 
primitive multiple ring of multiplicity $n$ associated to $K$).

\sepsub

\Ssect{Canonical filtrations}{mod_filt}

Let $M$ be a \m{R[n]}-module. Then it is easy to see that we have a canonical 
isomorphism \ \m{S_n^{-1}M\simeq M\ot_{R[n]}K[n]}, and a natural morphism of 
\m{R[n]}-modules \ \m{i_M:M\to M\ot_{R[n]}K[n]}.

For \m{0\leq i\leq n}, let \ \m{M_i=t^iM}, which is a submodule of $M$. The 
filtration
\[0=M_n\subset M_{n-1}\subset\cdots\subset M_1\subset M_0=M\]
is called the {\em first canonical filtration} of $M$. For \m{0\leq i<n}, let
\[G_i(M) \ = \ M_i/M_{i+1} \ \simeq \ M_i\ot_{R[n]}R \ . \]
For \m{0\leq i\leq n}, let \ \m{M^{(i)}=\{m\in M;t^im=0\}}, which is a 
submodule of $M$ that contains \m{M_{n-i}}. The filtration
\[M^{(0)}=\nsp\subset M^{(1)}\subset\cdots\subset M^{(n-1)}\subset M^{(n)}=M\]
is called the {\em second canonical filtration} of $M$. For \m{0<i\leq n}, let
\[G^{(i)}(M) \ = \ M^{(i)}/M^{(i-1)} \ , \]
which is a $R$-module.

For example we have, for \ \m{M=R[n]^p} \ and \m{0\leq i<n},
\begin{equation}\label{equ1}
M_i \ = \ M^{(n-i)} \ \simeq \ (R[n-i])^p \ , \qquad G_i(M) \ \simeq \ R^p \ . 
\end{equation}

For \m{1\leq i\leq n-1}, the multiplication by $t$, \m{M\to M}, induces an 
injective morphism of $R$-modules
\[\lambda_i:G^{(i+1)}(M)\lra G^{(i)}(M)\]
(for a $A$-module $M$ : \m{\lambda_i:G^{(i+1)}(M)\ot_A\ki\to G^{(i)}(M)}).
Let \ \m{\Gamma^{(i-1)}(M)=\coker(\lambda_i)}. Similarly we have for \ \m{0\leq 
i\leq n-2} \ a surjective morphism of $R$-modules
\[\mu_i:G_i(M)\lra G_{i+1}(M)\]
(resp. \m{\mu_i:G_i(M)\ot_A\ki\lra G_{i+1}(M)}). Let \ 
\m{\Gamma_i(M)=\ker(\mu_i)}.

We have then a canonical isomorphism \ \m{\Gamma_i(M)\simeq\Gamma^{(i)}(M)} \
(resp. \m{\Gamma_i(M)\simeq\Gamma^{(i)}(M)\ot_A\ki^{i+1}}, cf. prop. 
\ref{prop15}).

\end{sub}

\sepsub

\Ssect{Some canonical extensions}{canex}

We can view \m{R[i]}, \m{0\leq i<n}, as a \m{R[n]}-module. We have an exact 
sequence
\[0\lra R[i]=tR[i+1]\hookrightarrow R[i+1]\lra R\lra 0 \]
(resp. \ \m{0\to\ki A_{i+1}\simeq A_i\to A_{i+1}\to A\to 0}).

\sepprop

\begin{subsub}\label{lem2}{\bf Lemma: } We have \ 
\m{\Ext^1_{R[n]}(R,R[i])\simeq 
R}, and given an extension \Nligne \m{0\to R[i]\to P\to R\to 0} of 
\m{R[n]}-modules, corresponding to \m{\sigma\in R}, then \m{P\simeq R[i+1]} if 
and only if $\sigma$ is invertible.
\end{subsub}
\begin{proof}
We use \ref{extmod} and the following free resolution of  $R$ as a 
\m{R[n]}-module:
\xmat{\cdots R[n]\ar[rr]^-{Xt^{n-1}} & & R[n]\ar[rr]^-{Xt} & & R[n]\flon[r] & R 
\ 
.}
Let \ \m{\iota:R[n-1]=tR[n]\to R[n]} \ be the inclusion, and \ 
\m{\eta=\sigma\oplus\iota:R[n-1]\to R[i]\oplus R[n]}. Then \ 
\m{P=\coker(\eta)}. Lemma \ref{lem2} follows easily from this description.
\end{proof}

\end{sub}

\sepsub

\Ssect{Quasi free modules}{QFM}

For \m{1\leq i\leq n}, \m{R[i]} is a \m{R[n]}-module. Let $M$ be a 
\m{R[n]}-module.
We say that $M$ is {\em quasi free} if there exist integers \m{n_i\geq 0}, 
\m{1\leq i\leq n}, such that \ \m{M\simeq\bigoplus_{i=1}^nR[i]^{n_i}}.

The proof of the following result is similar to that of \cite{dr2}, 
th\'eor\`eme 5.1.3.

\sepprop

\begin{subsub}\label{theo1}{\bf Theorem: } Let $M$ be a \m{R[n]}-module of 
finite type. Then $M$ is quasi free if and only \m{G_i(M)} is a free $R$-module 
for \m{0\leq i\leq n}.
\end{subsub}

\sepprop

If $R$ is a field then the \m{G_i(M)} are $R$-vector spaces, hence they are 
free. It follows that

\sepprop

\begin{subsub}\label{cor1}{\bf Corollary: } Every finitely generated 
\m{K[n]}-module is quasi free (as a \m{K[n]}-module).
\end{subsub}

\end{sub}

\sepsub

\Ssect{Torsion}{tors}

Let $M$ be a \m{R[n]}-module of finite type.
The set of elements \m{m\in M} such that there exists \m{\alpha\in S_n} such 
that \m{\alpha m=0} is a submodule \m{T(M)} of $M$, called the {\em torsion 
submodule} of $M$. We say that $M$ is a {\em torsion module} if \ 
\m{M=T(M)}, and that $M$ is {\em torsion free} if \m{T(M)=0}, or 
equivalently if for every non zero \m{m\in M} and \m{u\in R[n]}, we have \ 
\m{um\not=0}. The module \m{M/T(M)} is torsion free. It is easy to see that \ 
\m{T(M)=\ker(i_M)} (cf. \ref{mod_filt}). Hence $M$ is torsion free if and only 
if \m{i_M} is injective.

\sepprop

\begin{subsub}\label{prop2}{\bf Proposition: } $M$ is torsion free if and only 
if it is isomorphic to a submodule of a free \m{R[n]}-module.
\end{subsub}
\begin{proof}
A free module is torsion free, so is $M$ if it is a submodule of a free module. 
Conversely, suppose that $M$ is torsion free. Then \ \m{i_M:M\to 
M\ot_{R[n]}K[n]} 
\ is injective. By corollary \ref{cor1}, \m{M\ot_{R[n]}K[n]} is quasi  
free as a \m{K[n]}-module, i.e. there are integers \m{n_i\geq 0}, \m{1\leq 
i\leq n}, such that we have an isomorphism
\[M\ot_{R[n]}K[n] \ \simeq \ \bigoplus_{i=1}^nK[i]^{n_i} \ . \]
Let \m{m_1,\ldots,m_q} be generators of $M$. We can write \
\m{i_M(m_j)=\sigg_{i=1}^n\lambda_{j,i}} \ , with \ \m{\lambda_{j,i}\in 
K[i]^{n_i}}. There exists \m{\mu\in R[n]}, such that \m{\mu_0\not=0}, and \ 
\m{\mu\lambda_{j,i}\in R[i]^{n_i}} \ for \m{1\leq j\leq q}, \m{1\leq i\leq n}. 
We have then
\[i_M(M) \ \subset \ \bigoplus_{i=1}^n\bigg(\frac{1}{\mu}R[i]\bigg)^{n_i} \ , \]
and since \m{\dsp\bigg(\frac{1}{\mu}R[i]\bigg)^{n_i}\simeq R[i]^{n_i}}, we 
obtain 
an inclusion \ \m{M\subset\bigoplus_{i=1}^nR[i]^{n_i}}. But \Nligne 
\m{R[i]\simeq (t^{n-i})\subset R[n]}, and finally \ \m{M\subset R[n]^m}, with 
\m{m=\sigg_{i=1}^nn_i}.
\end{proof}

\sepprop

\begin{subsub}\label{prop3}{\bf Proposition: } \m{\Ext^1_{R[n]}(M,R[n])} is a 
torsion module.
\end{subsub}
\begin{proof}
Using the exact sequence \ \m{0\to T(M)\to M\to M/T(M)\to 0}, it suffices to 
treat the cases of torsion modules, and of torsion free modules.

Suppose that \ \m{M=T(M)}. Since \m{T(M)} is finitely generated, there exists \
\m{x=\sigg_{i=0}^{n-1}u_it^i\in R[n]} \ with \m{u_0\not=0}, such that \ 
\m{xT(M)=\nsp}. The morphism \ 
\m{\Ext^1_{R[n]}(M,R[n])\to\Ext^1_{R[n]}(M,R[n])} \ 
induced by the multiplication by \m{x:M\to M}, is also the multiplication by 
$x$, and it is zero, i.e. \m{x\Ext^1_{R[n]}(M,R[n])=\nsp}.

Suppose that $M$ is torsion free. Let \ \m{\sigma\in\Ext^1_{R[n]}(M,R[n])}, 
corresponding to an extension
\xmat{0\ar[r] & R[n]\ar[r]^-\alpha & N\ar[r]^-p & M\ar[r] & 0 \ .}
The module $N$ is finitely generated and torsion free, hence by proposition 
\ref{prop2}, it is a submodule of a free module: \m{N\subset R[n]^m}. Let \ 
\m{\alpha(1)=(r_1,\ldots,r_m)\in R[n]^m}. Since $\alpha$ is injective, some 
\m{r_i} is not a multiple of $t$. Let \ \m{p_i:N\to R[n]} \ be the restriction 
to $N$ of the $i$th projection \m{R[n]^m\to R[n]}, so that \ 
\m{p_i\circ\alpha:R[n]\to R[n]} \ is the multiplication by \m{r_i}. We have 
then 
a commutative diagram
\xmat{0\ar[r] & R[n]\ar[r]^-\alpha\ar[d]^{\times r_i} & 
N\ar[r]\ar[d]^{p_i\oplus 
p} & M\ar[r]\fleq[d] & 0\\
0\ar[r] & R[n]\ar[r] & R[n]\oplus M\ar[r] & M\ar[r] & 0}
where the lower sequence is the trivial exact sequence. It is associated to 
\m{r_i\sigma}, hence \m{r_i\sigma=0}.
\end{proof}

\sepprop

\begin{subsub}\label{prop4}{\bf Proposition: } $M$ is torsion free if and only 
\m{M^{(1)}} is. In this case \m{G^{(i)}(M)} and \m{M/M^{(i)}} are torsion 
free for \m{1\leq i\leq n}.
\end{subsub}
\begin{proof}
If $M$ is torsion free then so is \m{M^{(1)}\subset M}. It $M$ is not torsion 
free, let \m{m\in M}, \m{m\not=0}, be such that there exists \m{\alpha\in 
R[n]}, 
not a zero divisor, such that \m{\alpha m=0}. Let $k$ be the biggest integer 
such that \m{z^km\not=0}. Then \m{z^km\in M^{(1)}} and \m{\alpha z^km=0}. Hence 
\m{M^{(1)}} is not torsion free.

Suppose that $M$ is torsion free. Let \m{\ov{m}\in M/M^{(1)}}, $\alpha$ a non  
zero divisor in \m{R[n]}, be such that \m{\alpha\ov{m}=0}. If \m{m\in M} is 
over $\ov{m}$, we have \ \m{\alpha m\in M^{(1)}}, i.e. \m{z\alpha m=0}. Since 
$M$ is torsion free, we have \m{zm=0}, \m{m\in M^{(1)}} and \m{\ov{m}=0}. This 
proves that \m{M/M^{(1)}} is torsion free.

The fact that \m{G^{(i)}(M)} and \m{M/M^{(i)}} are torsion free for \m{1\leq 
i\leq n} is easily proved by induction on $i$.
\end{proof}

\end{sub}

\sepsub

\Ssect{Duality}{dual}

Let $M$ be a \m{R[n]}-module of finite type, and \ \m{M^{\vee}=\Hom(M,R[n])} \ 
the {\em dual} of $M$. If \m{1\leq i\leq n} and $M$ is a \m{R[i]}-module, we 
have
\[M^\vee \ = \ \Hom(M,(t^{n-i})) \ = \ \Hom(M,R[i]) \ , \]
So the dual of $M$ as a \m{R[i]}-module is also \m{M^\vee}. In particular \ 
\m{(R[i])^\vee\simeq R[i]}. The dual of a \m{R[n]}-module is torsion free. If 
$M$ 
is a torsion module, then \m{M^\vee=\nsp}.

The following is immediate

\sepprop

\begin{subsub}\label{lem3}{\bf Lemma: } For \m{1\leq i<n}, we have \ 
\m{(M^\vee)^{(i)}\simeq(M/M_i)^\vee}.
\end{subsub}

\sepprop

Let \m{t_M:M\to M^{\vee\vee}} be the canonical morphism. If \m{m\in M}, 
\m{t_M(m)} is the linear form
\[\xymatrix@R=5pt{M^\vee\ar[r] & R[n]\\ \phi\fmaps[r] & \phi(m)} \ .\]
If $M$ is a free module, \m{t_M} is an isomorphism.

\sepprop

\begin{subsub}{\bf Proposition: }\label{prop5} We have \ \m{\ker(t_M)=T(M)}.
\end{subsub}
\begin{proof}
First we prove that \m{t_M} is injective if $M$ is torsion free.
By proposition \ref{prop2}, $M$ is then a submodule of a free module $E$. Let 
\m{m\in M}, \m{m\not=0}. Then there exists a linear form \m{\psi:E\to R[n]} 
such 
that \ \m{\psi(m)\not=0}. If \m{\phi\in M^\vee} is the restriction of $\psi$, 
we have \ \m{\phi(m)\not=0}.

Now for a general $M$, the exact sequence \ \m{0\to T(M)\to M\to M/T(M)\to 0} \ 
implies that \ \m{(M/T(M))^\vee\simeq M^\vee}. We have a commutative diagram
\xmat{0\ar[r] & T(M)\ar[r] & M\ar[r]\ar[d]^{t_M} & 
M/T(M)\ar[r]\ar[d]^{t_{M/T(M)}} & 0 \\
& & M^{\vee\vee}\ar[r]^-\simeq & (M/T(M))^{\vee\vee}}
and the result follows from the injectivity of \m{t_{M/T(M)}}.
\end{proof}

\sepprop

\begin{subsub}\label{prop6}{\bf Proposition: } \m{\coker(t_M)} is a torsion 
module.
\end{subsub}
\begin{proof} We first prove the following result: suppose that we have an 
exact sequence \Nligne \m{0\to N\to M\hfl{\pi}{} Q\to 0} \ of \m{R[n]}-modules, 
and that \m{\coker(t_N)} and \m{\coker(t_Q)} are torsion modules. We show that 
\m{\coker(t_M)} is a torsion module. We have an exact sequence
\[0\to Q^\vee\hfl{^t\pi}{} M^\vee\to N^\vee\hfl{\delta}{}\Ext^1_{R[n]}(Q,R[n]) 
\ 
. \]
By proposition \ref{prop3}, \m{\Ext^1_{R[n]}(Q,R[n])} is a torsion module. Let 
\ 
\m{U=\ker(\delta)}, \m{T=\imm(\delta)}. We have exact sequences
\[0\to Q^\vee\hfl{^t\pi}{} M^\vee\to U\to 0 \ , \quad 0\to U\to N^\vee\to T\to 
0 \ , \]
whence exact sequences
\[0\to U^\vee\to M^{\vee\vee}\hfl{^{tt}\pi}{}Q^{\vee\vee}\to T_2\to 0 \ ,
\quad 0\to N^{\vee\vee}\hfl{\mu}{} U^\vee\to T_3\to 0 \ , \]
for some torsion modules \m{T_2}, \m{T_3}. Let \ \m{V=\imm(^{tt}\pi)}. It 
follows that we have a commutative diagram with exact rows:
\xmat{0\ar[r] & N\ar[d]^\alpha\ar[r] & M\ar[d]^-{t_M}\ar[r] & 
Q\ar[d]^\beta\ar[r] & 0\\
0\ar[r] & U^\vee\ar[r] & M^{\vee\vee}\ar[r] & V\ar[r] & 0 }
hence an exact sequence
\[\coker(\alpha)\lra\coker(t_M)\lra\coker(\beta)\lra 0 \ . \]
Now $\alpha$ is the composition \ \m{N\hfl{t_N}{}N^{\vee\vee}\hfl{\mu}{}U^\vee},
so \m{\coker(\alpha)} is a torsion module. and the composition \ 
\m{Q\hfl{\beta}{} V\subset Q^{\vee\vee}} \ is \m{t_Q}, so \m{\coker(\beta)} is 
a torsion module. It follows that \m{\coker(t_M)} is also a torsion module.

Proposition \ref{prop6} is then easily proved by induction on the minimal 
number of generators of $M$ (or on $n$, using the canonical filtrations of $M$).
\end{proof}

\end{sub}

\sepsub

\Ssect{Balanced modules}{bal_mod}

Let $M$ be a \m{R[n]}-module of finite type. Let
\[\boldsymbol{\lambda}(M) \ = \ 
\lambda_1\circ\cdots\circ\lambda_{n-1}:G^{(n)}(M)\lra G^{(1)}(M) \ , \]
\[\boldsymbol{\mu}(M) \ = \ \mu_{n-2}\circ\cdots\circ\mu_0:G_0(M)\lra 
G_{n-1}(M)\]
(resp.
\[\boldsymbol{\lambda}(M) \ = \ 
\lambda_1\circ\cdots\circ\lambda_{n-1}:G^{(n)}(M)\ot_{A_n}\ki^{n-1}\lra 
G^{(1)}(M) \ , \]
\[\boldsymbol{\mu}(M) \ = \ 
\mu_{n-2}\circ\cdots\circ\mu_0:G_0(M)\ot_{A_n}\ki^{n-1}\lra G_{n-1}(M) \ ) \ . 
\]
(cf. \ref{mod_filt}).

\sepprop

\begin{subsub}\label{def2}{\bf Definition: } We say that $M$ is {\em balanced} 
if \m{\boldsymbol{\lambda}(M)} is surjective.
\end{subsub}

\sepprop

\begin{subsub}\label{prop18}{\bf Proposition:} The following properties are 
equivalent
\begin{enumerate}
\item[(i)] $M$ is balanced.
\item[(ii)] $\lambda_1,\ldots,\lambda_{n-1}$ are surjective.
\item[(iii)] $\Gamma_1(M)=\cdots=\Gamma_{n-1}(M)=0$ .
\item[(iv)] $\Gamma^{(1)}(M)=\cdots=\Gamma^{(n-1)}(M)=0$ .
\item[(v)] $\mu_0,\ldots,\mu_{n-2}$ are injective.
\item[(vi)] $\boldsymbol{\mu}(M)$ is injective.
\end{enumerate}
\end{subsub}
(cf. prop. \ref{prop16}).

\sepprop

\begin{subsub}\label{prop19}{\bf Proposition:} The \m{R[n]}-module $M$ is 
balanced if and only if \ \m{M_i=M^{(n-i)}} \ for \ \m{1\leq i\leq n}.
\end{subsub}
(cf. prop. \ref{prop17}).

\sepprop

In particular a \m{R[2]}-module is balanced if and only if \ \m{M^{(1)}=M_1}.

\sepprop

\begin{subsub}\label{prop30}{\bf Proposition: } Let $M$ be a balanced 
\m{R[n]}-module.

{\bf 1-- } Let
\[0=N_n\subset N_{n-1}\subset\cdots\subset N_1\subset N_0=N\]
be a filtration such that, for \m{0\leq i<n}, \m{N_i/N_{i+1}\not=\nsp} \ and \ 
\m{t.(N_i/N_{i+1})=\nsp}. Then we have \ \m{N_i=M_i} \ for \m{1\leq i<n}.

{\bf 2-- } Let \m{M'\subset M} be a submodule. Suppose that \ 
\m{t^{n-1}.M'=\nsp} \ and \ \m{t.(M/M')=\nsp}. Then we have \m{M'=M_1}.
\end{subsub}
\begin{proof} We first prove {\bf 1}. Let \m{u\in N_1}. Let \ \m{p:N_1
\to 
N_1/N_2} \ be the projection. Then \ \m{p(tu)=tp(u)=0}, hence \m{tu\in N_2}, 
and similarly \ \m{t^{i-1}u\in N_i} \ for \m{2\leq i<n}. It follows that \ 
\m{N_1\subset M^{(n-1)}}. Since \ \m{M_1\subset N_1} \ and \ \m{M_1=M^{(n-1)}}, 
we have \ \m{N_1=M_1}. Then {\bf 1} is easily proved by induction on $n$.

{\bf 2} is a consequence of {\bf 1} (take the first canonical filtration of 
\m{M'}).
\end{proof}

\end{sub}

\sepsub

\Ssect{Balanced ideals and regular sequences}{reg_seq}

\begin{subsub}\label{prop20}{\bf Proposition: } Let $p$ be a positive integer, 
and for \m{1\leq i\leq p}, \m{x_i=\sigg_{j=0}^{n-1}x_{i,j}t^j\in R[n]}, with \ 
\m{x_{i,j}\in R} \ for \m{0\leq j<n}. Let $I$ (resp. \m{I_0}) be the ideal of 
\m{R[n]} (resp. $R$) generated by \m{x_1,\ldots,x_p} (resp. 
\m{x_{1,0},\ldots,x_{p,0}}). Then
\begin{enumerate}
\item[(i)] $(x_1,\ldots,x_p)$ is a regular sequence in $R[n]$ if and only if 
$(x_{1,0},\ldots,x_{p,0})$ is a regular sequence in $R$.
\item[(ii)] If $(x_1,\ldots,x_p)$ is a regular sequence in $R[n]$, then for 
every $y\in R$, $t^{n-1}y\in I$ if and only if $y\in I_0$.
\end{enumerate}
\end{subsub}
\begin{proof}
By induction on $n$. The case \m{n=1} is obvious. Suppose that the theorem is 
true for \m{n-1\geq 1}. We will prove that it is true for $n$ by induction on 
$p$, the case \m{p=1} being obvious. Suppose that it is true for \m{p-1\geq 1}.

Suppose that \m{(x_1,\ldots,x_p)} is a regular sequence in \m{R[n]}. We will 
prove that \m{(x_{1,0},\ldots,x_{p,0})} is a regular sequence in $R$. By the 
induction hypothesis, \m{(x_{1,0},\ldots,x_{p-1,0})} is a regular sequence in 
$R$, and we have to show that \m{x_{p,0}} is not a zero divisor in 
\m{R/(x_{1,0},\ldots,x_{p-1,0})}. Let \m{a\in R} be such that \ 
\m{ax_{p,0}\in(x_{1,0},\ldots,x_{p-1,0})}. Hence \ 
\m{t^{n-1}ax_p\in(x_1,\ldots,x_{p-1})}. By (ii) and the induction hypothesis, 
we have \ \m{a\in(x_{1,0},\ldots,x_{p-1,0})}, i.e. \m{a=0} in \ 
\m{R/(x_{1,0},\ldots,x_{p-1,0})}. Hence \m{(x_{1,0},\ldots,x_{p,0})} is a 
regular sequence in $R$.

Suppose now that \m{(x_{1,0},\ldots,x_{p,0})} is a regular sequence in $R$. We 
will prove that \m{(x_1,\ldots,x_p)} is a regular sequence in \m{R[n]}. By the 
induction hypothesis, \m{(x_1,\ldots,x_{p-1})} is a regular sequence in 
\m{R[n]}, and we must show that \m{x_p} is not a zero divisor in 
\m{R[n]/(x_1,\ldots,x_{p-1})}. Let \m{a\in R[n]} be such that \ 
\m{ax_p\in(x_1,\ldots,x_{p-1})}. By the induction hypothesis (the result is 
true for \m{n-1}), we can write \ \m{a=b+t^{n-1}c}, with \ 
\m{b\in(x_1,\ldots,x_{p-1})} \ and \m{c\in R}. We have then \ 
\m{t^{n-1}cx_p\in(x_1,\ldots,x_{p-1})}. By (ii) and the induction hypothesis, 
we have \ \m{cx_{p,0}\in(x_{1,0},\ldots,x_{p-1,0})}. Since 
\m{(x_{1,0},\ldots,x_{p-1,0})} is a regular sequence, we have 
\m{c\in(x_{1,0},\ldots,x_{p-1,0})}. It follows that \ 
\m{a\in(x_1,\ldots,x_{p-1})}, hence \m{(x_1,\ldots,x_p)} is a regular sequence 
in \m{R[n]}.

Now we prove (ii). If \m{y\in I_0}, it is clear that \ \m{t^{n-1}y\in I}. 
Conversely, suppose that \ \m{t^{n-1}y\in I}. We can write 
\m{t^{n-1}y=t^k\sigg_{i=1}^p\alpha_ix_i}, with \m{k\in\{0,\ldots,n-1\}} 
maximal, and \ \m{\alpha_1,\ldots,\alpha_p\in R[n]}. We must show that 
\m{k=n-1}.
Suppose that \ \m{k<n-1}. We can write \ 
\m{t^{n-1-k}y=\sigg_{i=1}^p\alpha_ix_i+t^{n-k}z}, for some \m{z\in R[n]}, i.e. \
\m{t^{n-1-k}(y-tz)=\sigg_{i=1}^p\alpha_ix_i}. By the induction hypothesis (the 
result is true for \m{n-k}), the term of degree 0 of \m{y-tz} belongs to 
\m{I_0}, i.e \m{y\in I_0}.
\end{proof}

\sepprop

\begin{subsub}\label{prop21}{\bf Proposition: } Let \m{(x_1,\ldots,x_p)} be a 
regular sequence in \m{R[n]}, and $I$ the ideal of \m{R[n]} generated by 
\m{x_1,\ldots,x_p}. Then $I$ is a balanced \m{R[n]}-module.
\end{subsub}
\begin{proof}
Let \m{I'\subset R[n-1]} be the ideal generated by the images 
\m{x'_1,\ldots,x'_p} of \m{x_1,\ldots,x_p} in \m{R[n-1]}. By proposition 
\ref{prop20} (i), \m{(x'_1,\ldots,x'_p)} is a regular sequence in \m{R[n-1]}.
Let \ \m{\phi:I\to I'} \ be the restriction of the canonical morphism 
\m{R[n]\to R[n-1]}. Then we have \ \m{\ker(\phi)=I\cap(t^{n-1})}, and by 
proposition \ref{prop20} (ii), \m{I\cap(t^{n-1})=t^{n-1}I}. It follows 
that $\phi$ induces an isomorphism \ \m{I/t^{n-1}I\simeq I'}. Hence, using an 
induction on $n$, we see that for \m{1\leq k<n}, \m{I_k} is canonically 
isomorphic to the ideal of \m{R[k]} generated by the images of 
\m{x_1,\ldots,x_p}. Proposition \ref{prop21} follows easily.
\end{proof}

\sepprop

\begin{subsub}\label{ot-ex} Other examples -- \rm Suppose that \ \m{R=\C[X,Y]}, 
and let \m{I=(X^2,Y^2,XY)},\Nligne \m{J=(X^2,Y^2+t,XY)} in \m{R[2]}. Then we 
have \ \m{I_1=I^{(1)}=t.(X^2,Y^2,XY)}, so $I$ is balanced. But \ 
\m{J_1=t.(X^2,Y^2,XY)}, and \ \m{tX\in J^{(1)}\backslash J_1} 
(because \m{tX=X(Y^2+t)-Y.XY}). Hence $J$ is not balanced.
\end{subsub}

\end{sub}

\sepsec

\section{Primitive multiple schemes}\label{PMS}

\Ssect{Definition and construction}{PMS_def}

Let $X$ be a smooth connected variety, and \ \m{d=\dim(X)}. A {\em multiple 
scheme with support $X$} is a Cohen-Macaulay scheme $Y$ such that 
\m{Y_{red}=X}. If $Y$ is quasi-projective we say that it is a {\em multiple 
variety with support $X$}. In this case $Y$ is projective if $X$ is.

Let $n$ be the smallest integer such that \m{Y=X^{(n-1)}}, \m{X^{(k-1)}}
being the $k$-th infinitesimal neighborhood of $X$, i.e. \
\m{\ki_{X^{(k-1)}}=\ki_X^{k}} . We have a filtration \ \m{X=X_1\subset
X_2\subset\cdots\subset X_{n}=Y} \ where $X_i$ is the biggest Cohen-Macaulay
subscheme contained in \m{Y\cap X^{(i-1)}}. We call $n$ the {\em multiplicity}
of $Y$.

We say that $Y$ is {\em primitive} if, for every closed point $x$ of $X$,
there exists a smooth variety $S$ of dimension \m{d+1}, containing a 
neighborhood of $x$ in $Y$ as a locally closed subvariety. In this case, 
\m{L=\ki_X/\ki_{X_2}} is a line bundle on $X$, \m{X_j} is a primitive multiple 
scheme of multiplicity $j$ and we have \ 
\m{\ki_{X_j}=\ki_X^j}, \m{\ki_{X_{j}}/\ki_{X_{j+1}}=L^j} \ for \m{1\leq j<n}. 
We call $L$ the line bundle on $X$ {\em associated} to $Y$. The ideal sheaf 
\m{\ki_{X,Y}} can be viewed as a line bundle on \m{X_{n-1}}.

Let \m{P\in X}. 
Then there exist elements \m{y_1,\ldots,y_d}, $t$ of \m{m_{S,P}} whose images 
in \m{m_{S,P}/m_{S,P}^2} form a basis, and such that for \m{1\leq i<n} we have 
\ \m{\ki_{X_i,P}=(t^{i})}. In this case the images of \m{y_1,\ldots,y_d} in 
\m{m_{X,P}/m_{X,P}^2} form a basis of this vector space.

Even if $X$ is projective, we do not assume that $Y$ is projective.

The simplest case is when $Y$ is contained in a smooth variety $S$ of dimension 
\m{d+1}. Suppose that $Y$ has multiplicity $n$. Let \m{P\in X} and 
\m{f\in\ko_{S,P}}  a local equation of $X$. Then we have \ 
\m{\ki_{X_i,P}=(f^{i})} \ for \m{1<i\leq n} in $S$, in particular 
\m{\ki_{Y,P}=(f^n)}, and \ \m{L=\ko_X(-X)} .

For any \m{L\in\Pic(X)}, the {\em trivial primitive variety} of multiplicity 
$n$, with induced smooth variety $X$ and associated line bundle $L$ on $X$ is 
the $n$-th infinitesimal neighborhood of $X$, embedded by the zero section in 
the dual bundle $L^*$, seen as a smooth variety.

\sepsubsub

\begin{subsub}\label{PMS-1} Construction of primitive multiple schemes -- \rm
Let $Y$ be a primitive multiple scheme of multiplicity $n$, \m{X=Y_{red}}.
Let \ \m{{\bf Z}_n=\spec(\C[t]/(t^n))}.
Then for every closed point \m{P\in X}, there exists an open neighborhood $U$ 
of $P$ in $X$, and an open neighborhood \m{U^{(n)}} of $P$ in $Y$ such that 
\begin{enumerate}
\item[--] $U^{(n)}\cap X=U$ ,
\item[--] There exists a section \ \m{\ko_X(U)\to\ko_Y(U^{(n)})} \ of the 
restriction map \ \m{\ko_Y(U^{(n)})\to\ko_X(U)},
\item[--] $L_{|U}$ is trivial,
\end{enumerate}
and there exists a commutative diagram
 \xmat{ & U\flinc[ld]\flinc[rd] \\
U^{(n)}\ar[rr]^-\simeq & & U\times {\mathbf Z}_n}
i.e. $Y$ is locally trivial (\cite{dr1}, th\'eor\`eme 5.2.1, corollaire 5.2.2).
For every open subset $V$ of $X$, \m{V^{(n)}} will denote the corresponding 
open subset of $Y$.

It follows that we can construct a primitive multiple scheme of multiplicity 
$n$ by taking an open cover \m{(U_i)_{i\in I}} of $X$ and gluing the varieties 
\ \m{U_i\times{\bf Z}_n} (with automorphisms of the \ \m{U_{ij}\times{\bf Z}_n} 
\ leaving \m{U_{ij}} invariant).

Let \m{(U_i)_{i\in I}} be an affine open cover of $X$ such that we have 
trivializations
\xmat{\delta_i:U_i^{(n)}\ar[r]^-\simeq & U_i\times{\bf Z}_n , }
and \ \m{\delta_i^*:\ko_{U_i\times{\bf Z}_n}\to\ko_{U_i^{(n)}}} \ the 
corresponding isomorphisms. Let
\xmat{\delta_{ij}=\delta_j\delta_i^{-1}:U_{ij}\times{\bf Z}_n\ar[r]^-\simeq & 
U_{ij}\times{\bf Z}_n \ . }
Then \ \m{\delta_{ij}^*=\delta_i^{*-1}\delta_j^*} \ is an automorphism of \ 
\m{\ko_{U_i\times Z_n}=\ko_X(U_{ij})[t]/(t^n)}, such that for every \ 
\m{\phi\in\ko_X(U_{ij})}, seen as a polynomial in $t$ with coefficients in 
\m{\ko_X(U_{ij})}, the term of degree zero of \m{\delta_{ij}^*(\phi)} is the 
same as the term of degree zero of $\phi$.
\end{subsub}

\sepsubsub

\begin{subsub}\label{I_X} The ideal sheaf of $X$ -- \rm
There exists \ \m{\alpha_{ij}\in\ko_X(U_{ij})\ot_\C(\C[t]/(t^{n-1}))} \ such 
that \ \m{\delta_{ij}^*(t)=\alpha_{ij}t}. Let \ 
\m{\alpha^{(0)}_{ij}=\alpha_{ij|X}\in\ko_X(U_i)}.
For every \m{i\in I}, \m{\delta_i^*(t)} is a generator of 
\m{\ki_{X,Y|{U^{(n)}}}}. So we have local trivializations
\[\xymatrix@R=5pt{\lambda_i:\ki_{X,Y|{U_i^{(n-1)}}}\ar[r] & 
\ko_{U_i^{(n-1)}}\\
\delta_i^*(t)\fmaps[r] & 1}\]
Hence \ \m{\lambda_{ij}=\lambda_i\lambda_j^{-1}: 
\ko_{U_{ij}^{(n-1)}}\to\ko_{U_{ij}^{(n-1)}}} \ is the multiplication by 
\m{\delta_j^*(\alpha_{ij})}. It follows that 
\m{(\delta_j^*(\alpha_{ij}))} (resp. 
\m{(\alpha^{(0)}_{ij})})  is a cocycle representing the line bundle 
\m{\ki_{X,Y}} (resp. \m{L}) on \m{X_{n-1}} (resp. $X$). 
\end{subsub}

\end{sub}

\sepsub

\Ssect{The case of double schemes}{DBC}

We suppose that \m{n=2}. Let \ \m{\alpha_i:L_{|U_i}\to\ko_{U_i}} \ be 
isomorphisms such that \ \m{\alpha_{ij}=\alpha_i\circ\alpha_j^{-1}} \ on 
\m{U_{ij}}. Then we have the following description of \m{\delta_{ij}^*}:
\begin{enumerate}
\item[--] There are derivations 
\m{D_{ij}} of \m{\ko_X(U_{ij})} such that \ \m{\delta_{ij}^*(\beta)= 
\beta+D_{ij}(\beta)t} \ for every \m{\beta\in\ko_X(U_{ij})}.
\item[--] $\delta_{ij}^*(t)=\alpha_{ij}t$ .
\end{enumerate}
The relation \ \m{\delta_{ij}^*\delta_{jk}^*=\delta_{ik}^*} \ is equivalent to 
\m{D_{ij}+\alpha_{ij}D_{jk}=D_{ik}}. We can view \m{D_{ij}} as section of 
\m{T_{X|U_{ij}}}, the family \m{(\alpha_i^{-1}\ot D_{ij})} represents an 
element $\lambda$ of \m{H^1(X,T_X\ot L)}, and \m{\C\lambda} is independent of 
the choice of the \m{\delta_{ij}^*} and \m{\alpha_i}. We will note \ 
\m{\C\lambda=\zeta(X_2)}.

According to \cite{dr10} (and \cite{ba_ei}), two primitive double 
schemes \m{X_2}, \m{X'_2}, with underlying smooth variety $X$ and associated 
line bundle $L$ are isomorphic (over $X$) if and only if 
\m{\zeta(X_2)=\zeta(X'_2)}. And \m{X_2} is the trivial primitive double 
scheme 
if and only \ \m{\zeta(X_2)=0}.

\end{sub}

\sepsub

\Ssect{Canonical filtrations, quasi locally free sheaves}{QF}

Let $X$ be a smooth and irreducible variety. Let \m{Y=X_n} be a
primitive multiple scheme of multiplicity $n$, with underlying smooth variety
$X$, and associated line bundle $L$ on $X$.

Let \m{P\in X} be a closed point.

\sepprop

The two {\em canonical filtrations} are useful tools to study the coherent
sheaves on primitive multiple schemes. They have been defined for 
\m{\ko_{X_n,P}}-modules in \ref{mod_filt}.

\sepprop

\begin{subsub}\label{QLL-def1} \rm The {\em first canonical filtration of 
$\ke$} is
\[\ke_n=0\subset \ke_{n-1}\subset\cdots\subset \ke_{1}\subset \ke_0=\ke\]
where for \m{0\leq i< n}, \m{\ke_{i+1}} is the kernel of the canonical
surjective morphism \ \m{\ke_i\to\ke_{i\mid X}}. For \m{0\leq i<n}, let
\[G_i(\ke)\ = \ \ke_{i}/\ke_{i+1} \ = \ \ke_{i\mid X} \ . \]
We have \ \m{\ke/\ke_i=\ke_{\mid X_i}} . 
\end{subsub}

\sepsubsub

\begin{subsub}\label{2-fc} \rm One defines similarly the {\em second canonical 
filtration of $\ke$}:
\[\ke^{(0)}=\nsp\subset \ke^{(1)}\subset\cdots\subset
\ke^{(n-1)}\subset \ke^{(n)}=\ke  ,\]
where, for \m{1\leq i\leq n-1}, \m{\ke^{(i)}} is the maximal subsheaf $\kf$ of 
$\ke$ such that \ \m{\ki_{X,X_n}^i\kf=0}. For \m{0<i\leq n}, let
\[G^{(i)}(\ke) \ = \ \ke^{(i)}/\ke^{(i-1)} \ . \]
\end{subsub}

\sepsubsub

\begin{subsub}\label{prop_filt} Relationship between the two filtrations -- \rm
Let $\ke$ be a coherent sheaf on \m{X_n}. For \ \m{1\leq i\leq n-1}, the 
canonical morphism \ \m{\ke\ot\ki_{X,X_n}\to\ke} \ induces an injective 
morphism of sheaves on $X$
\[\lambda_i:G^{(i+1)}(\ke)\ot L\lra G^{(i)}(\ke) \ . \]
Let \ \m{\Gamma^{(i-1)}(\ke)=\coker(\lambda_i)}. Similarly we have for \ 
\m{0\leq i\leq n-2} \ a surjective morphism of sheaves on $X$
\[\mu_i:G_i(\ke)\ot L\lra G_{i+1}(\ke) \ . \]
Let \ \m{\Gamma_i(\ke)=\ker(\mu_i)}.
\end{subsub}

\sepprop

\begin{subsub}\label{prop15}{\bf Proposition: } For \ \m{0\leq i<n} \ we have \
\m{\Gamma_i(\ke)\simeq\Gamma^{(i)}(\ke)\ot L^{i+1}}.
\end{subsub}
\begin{proof}
We suppose that \m{\ki_{X,X_n}}, which is a line bundle on \m{X_{n-1}}, can be 
extended to a line bundle $\L$ on \m{X_n}. Let \ 
\m{\phi_i:\ke_i\ot\L\to\ke_{i+1}} \ be the canonical 
surjective morphism. Then we have a commutative diagram with exact rows
\xmat{0\ar[r] & \ke^{(i)}\ar[r]\flinc[d] & \ke\ar[r]\fleq[d] & 
\ke_i\ot\L^{-i}\ar[r]\flon[d]^{\phi_i\ot I_{\L^{-i-1}}} & 0\\
0\ar[r] & \ke^{(i+1)}\ar[r] & \ke\ar[r] & \ke_{i+1}\ot\L^{-i-1}\ar[r] & 0
}
It follows that \ \m{\ker(\phi_i)\simeq G^{(i+1)}(\ke)\ot L^{i+1}}.

The result follows from the commutative diagram with two exact rows and columns
\xmat{& & 0\ar[d] & 0\ar[d]\\
0\ar[r] & G^{(i+2)}(\ke)\ot L^{i+2}\ar[r]\flinc[d]^{\lambda_{i+1}\ot 
I_{L^{i+1}}} &
\ke_{i+1}\ot\L\ar[r]\ar[d] & \ke_{i+2}\ar[d]\ar[r] & 0\\
0\ar[r] & G^{(i+1)}(\ke)\ot L^{i+1}\ar[r] & \ke_i\ot\L\ar[r]\ar[d] &
\ke_{i+1}\ar[r]\ar[d] & 0\\
& & G_i(\ke)\ot L\ar[r]^-{\mu_i}\flon[r]\ar[d] & G_{i+1}(\ke)\ar[d]\\
& & 0 & 0}

In the general case, $X$ can be covered by open subsets $U$ such that, if  
\m{U_n} is the open subset of \m{X_n} corresponding to $U$, then
\m{\ki_{X,X_n|U_n}} can be extended to a line bundle on \m{U_n}. We obtain 
isomorphisms \ \m{\Gamma_i(\ke_{|U_n})\simeq\Gamma^{(i)}(\ke_{|U_n})\ot 
L^{i+1}}, and it is easy to see that they coincide on the intersections of any 
two of these open subsets $U$, and thus define the global isomorphism of 
proposition \ref{prop15}.
\end{proof}

\sepprop

\begin{subsub}\label{qlf} Quasi locally free sheaves -- \rm Let $P$ be a closed 
point of $X$.
Let $M$ be a \m{\ko_{X_n,P}}-module of finite type. Then $M$ is called {\em
quasi free} if there exist non negative integers \m{m_1,\ldots,m_n} and an
isomorphism \ \m{M\simeq\oplus_{i=1}^nm_i\ko_{X_i,P}}~. The integers
\m{m_1,\ldots,m_n} are uniquely determined: it is easy to recover them from 
the first canonical filtration of $M$. We say that \m{(m_1,\ldots,m_n)} is the 
{\em type} of $M$.

Let $\ke$ be a coherent sheaf on \m{X_n}. We say that $\ke$ is {\em quasi free
at $P$} if \m{\ke_P} is quasi free, and that $\ke$ is {\em quasi locally free} 
on a nonempty open subset \m{U\subset X} if it is quasi free at every point of 
$U$. If \m{U=X} we say that $\ke$ is quasi locally free. In this case the types 
of the modules \m{\ke_x}, \m{x\in X}, are the same (the {\em type} of $\ke$), 
and for every \m{x\in X} there exists a neighborhood \m{V\subset X_n} of $x$ 
and an isomorphism
\begin{equation}\label{equ20}\ke_{|V} \ \simeq \ \bigoplus _{n=1}^n 
m_i\ko_{X_i\cap 
V} \ .
\end{equation}
\end{subsub}

\sepprop

For every coherent sheaf $\ke$ on \m{X_n} there exists a nonempty open subset 
\m{U\subset X} such that \m{\ke} is quasi locally free on $U$ (cf. \cite{dr11}, 
4.7).

\sepprop

\begin{subsub}\label{prop10}{\bf Proposition: } Let $\ke$ be a quasi locally 
free sheaf on \m{X_n}. Then 

{\rm (i) } $\ke$ is reflexive. 

{\rm (ii)} For every positive integer 
$i$ we have \ \m{\EExt^i_{\ko_{X_n}}(\ke,\ko_{X_n})=0}.

{\rm (iii)} For every vector bundle $\F$ on \m{X_n} and every integer $i$ we 
have \Nligne \m{\Ext^i_{\ko_{X_n}}(\ke,\F)\simeq H^i(X_n,\ke^\vee\ot\F)} .

{\rm (iv)} Let \m{Y\subset X} be a closed subvariety of codimension \m{\geq 2}, 
\m{U=X\backslash Y}. Let $\ke$ be a coherent sheaf on \m{X_n}, such that 
\m{\ke_P} is quasi free for every \m{P\in Y}. Then the restriction map \ 
\m{H^0(X_n,\ke)\to H^0(U^{(n)},\ke)} \ is an isomorphism.
\end{subsub}

\begin{proof}
The assertion (iii) follows easily from (i), with the Ext spectral sequence,
and (iv) is immediate. 
Assertions (i) and (ii) are local, so we can replace \m{X_n} with \ \m{U\times 
{\bf Z}_n}, where $U$ is an open subset of \m{X_n}. Let \m{U_i=U\times{\bf 
Z}_i}, for \m{1\leq i\leq n}. The results follow from $(\ref{equ20})$ and the 
free resolutions
\xmat{....\ar[r] & \ko_{U_n}\ar[r]^-{\times t^{i}}\ar[r] & 
\ko_{U_n}\ar[r]^-{\times t^{n-i}}\ar[r] & \ko_{U_n}\ar[r]^-{\times t^{i}}
\ar[r] & \ko_{U_i}}
\end{proof}

\end{sub}

\sepsub

\Ssect{Morphisms of quasi locally free sheaves}{mor_QLL}

Let \m{r,r_1,\ldots,r_n} be integers, with \m{r_i\geq 0} and \ 
\m{r\geq\sum_{i=1}^nr_i}. Let \ 
\m{\ke=\bigoplus_{1\leq i\leq n}(\ko_{X_i}\ot\C^{r_i})} \ and \ 
\m{\phi:\ko_{X_n}\ot\C^r\to\ke} \ a morphism. 

Let \ \m{k=(\sum_{i=1}^nr_i)-r}. Chose a direct sum decomposition \ 
\m{\C^r=\C^k\oplus\big(\bigoplus_{1\leq i\leq n}\C^{r_i}\big)}. The {\em 
canonical morphism} \ \m{\phi_0:\ko_{X_n}\ot\C^r\to\ke} \ is defined as 
follows: \m{\phi_0=0} on \m{\ko_{X_n}\ot\C^k}, and on \m{\ko_{X_n}\ot\C^{r_i}},
\m{\phi=p_i\ot I_{\C^{r_i}}}, where \ \m{p_i:\ko_{X_n}\to\ko_{X_i}} is the 
restriction.

\sepprop

\begin{subsub}\label{lem20}{\bf Lemma: } {\bf 1 -- } Let \m{x\in X} be a closed 
point. If $\phi$ is surjective at $x$ then there exists an open neighborhood 
$U$ of $x$, and automorphisms $\alpha$ of \ \m{\ko_{U^{(n)}}\ot\C^r}, $\beta$ 
of \m{\ke_{|U^{(n)}}} such that \ 
\m{\phi_{|U^{(n)}}=\beta\circ\phi_{0|U^{(n)}}\circ\alpha}.

{\bf 2 -- } $\phi$ is surjective if and only of \m{\phi_{|X}} is. In this case 
\m{\ker(\phi)} is quasi locally free.
\end{subsub}
\begin{proof} We will prove {\bf 1} by induction on $n$, and {\bf 2} is an easy 
consequence of {\bf 1}. The result is obvious if \m{n=1}. Suppose that \m{n>1} 
and that the lemma is true on \m{X_{n-1}}. The component of $\phi$, 
\m{\phi_n:\ko_{X_n}\ot\C^r\to\ko_{X_n}\ot\C^{r_n}}, is surjective, so we have a 
direct sum decomposition \ \m{\C^r=\C^{r_n}\oplus\C^{r-r_n}} such that there 
exist $U$, $\alpha$, $\beta$ such that the component \m{\phi'_n} of \ 
\m{\phi'=\beta^{-1}\circ\phi_{|U^{(n)}}\circ\alpha^{-1}} \ is the identity on 
\m{\ko_{U^{(n)}}\ot\C^{r_n}} and zero on \m{\ko_{U^{(n)}}\ot\C^{r-r_n}}. We can 
even choose $\alpha$, $\beta$ such that \ \m{\phi'(\ko_{U^{(n)}}\ot\C^{r_n})} \ 
is the summand \ \m{\ko_{U^{(n)}}\ot\C^{r_n}} \ of \m{\ke_{|U^{(n)}}}. On the 
summand \ \m{\ko_{U^{(n)}}\ot\C^{r-r_n}} \ of \m{\ko_{X_n}\ot\C^r}, \m{\phi'} 
induces a morphism \Nligne 
\m{\phi'':\ko_{U^{(n)}}\ot\C^{r-r_n}\to\bigoplus_{1\leq 
i\leq n-1}(\ko_{X_i}\ot\C^{r_i})} that vanishes on \ 
\m{L^{n-1}\ot\C^{r-r_n}}, so that \m{\phi''} is a surjective morphism \ 
\m{\ko_{U^{(n-1)}}\ot\C^{r-r_n}\to\bigoplus_{1\leq i\leq 
n-1}(\ko_{X_i}\ot\C^{r_i})}. The result is easily obtained by applying the 
induction hypothesis on \m{\phi''}.
\end{proof}

\sepprop

\begin{subsub}\label{cor5}{\bf Corollary: } Let $\kn$, $\ke$ coherent sheaves 
on \m{X_n}, with $\kn$ quasi locally free. Let
\[0\lra\kn\hfl{\psi}{}\ko_{X_n}\ot\C^r\lra\ke\lra 0\]
an exact sequence. The $\ke$ is quasi locally free if and only if
\[{}^t\psi:\ko_{X_n}\ot\C^r\lra\kn^*\]
is surjective.
\end{subsub}
\begin{proof}
If $\ke$ is quasi locally free, then \ \m{\EExt^1_{\ko_{X_n}}(\ke,\ko_{X_n})=0} 
\ by proposition \ref{prop10}, hence \m{{}^t\psi} is surjective.

Conversely, suppose that \m{{}^t\psi} is surjective. We have then an exact 
sequence
\[0\lra\ke^*\lra\ko_{X_n}\ot\C^r\hfl{{}^t\psi}{}\kn^*\lra 0 \ , \]
hence \m{\ke^*} is quasi locally free by lemma \ref{lem20}. Since 
\m{\EExt^1_{\ko_{X_n}}(\kn^*,\ko_{X_n})=0} \ by proposition \ref{prop10},
we have an exact sequence
\[0\lra\kn^{**}=\kn\hfl{\psi}{}\ko_{X_n}\ot\C^r\lra\ke^{**}\lra 0 \ . \]
It follows that \ \m{\ke^{**}=\ke}, and $\ke$ is quasi locally free.
\end{proof}

\end{sub}

\sepsub

\Ssect{Flat families of quasi locally free sheaves}{fl_QLL}

Let $C$ be an irreducible smooth curve.

\sepprop

\begin{subsub}\label{theo4} {\bf Theorem: } Let $\ke$ be a coherent sheaf on 
\m{X_n\times C}, flat on $C$. Suppose that for every closed point \m{c\in C}, 
\m{\ke_c} is quasi locally free, and that the type of \m{\ke_c} is independent 
of $c$. Then 

{\bf 1 -- } $\ke$ is quasi locally free.

{\bf 2 -- } For \m{1\leq i\leq n-1}, \m{\ke_i} and \m{\ke_{i-1}/\ke_i} are flat 
on $C$, and for every \m{c\in C}, the morphism \ \m{(\ke_i)_c\to(\ke_{i-1})_c} 
\ induced by \ \m{\ke_i\subset\ke_{i-1}} is injective, and
\[(\ke_i)_c \ = \ (\ke_c)_i \ , \quad \ (\ke_{i-1}/\ke_i)_c \ = 
(\ke_{i-1})_c/(\ke_i)_c \ . \]

{\bf 3 -- } For \m{1\leq i\leq n-1}, \m{\ke^{(i)}} and 
\m{\ke^{(i+1)}/\ke^{(i)}} are flat on $C$, and for every \m{c\in C},the 
morphism \ \m{(\ke^{(i)})_c\to(\ke^{(i+1)})_c} \ induced by \ 
\m{\ke^{(i)}\subset\ke^{(i+1)}} is injective, and
\[(\ke^{(i)})_c \ = (\ke_c)^{(i)} \ , \quad (\ke^{(i+1)}/\ke^{(i)})_c \ =
(\ke^{(i+1)})_c/(\ke^{(i)})_c \ . \]
\end{subsub}
\begin{proof} We will prove the same statement with a nonempty subset of 
\m{X\times C} instead of \m{X\times C}. The proof is by induction on $n$. The 
result is obvious if \m{n=1}. Suppose that \m{n>1} and that it is true for 
\m{n-1}.

Let \m{x\in X}, \m{c_0\in C} be closed points. Since \m{\ke_{c_0}} is quasi 
locally free, there exists an open neighborhood \ \m{V\subset X} \ of $x$ such 
that \m{\ke_{c_0|V}} is free: \m{\ke_{c_0|V}\simeq \ko_V\ot\C^r}. It 
follows that there exists an open neighborhood \ \m{U\subset X\times C} \ of 
\m{(x,c_0)} and a surjective morphism
\[\phi:\ko_{U^{(n)}}\ot\C^r\lra\ke_{|U^{(n)}} \ . \]
Let \m{\kn=\ker(\phi)}. From the exact sequence
\[0\lra\kn\hfl{\psi}{}\ko_{U^{(n)}}\ot\C^r\hfl{\phi}{}\ke_{|U^{(n)}}\lra 0\]
and \cite{sga1}, Expos\'e IV, proposition 1.1, $\kn$ is flat on $C$. For 
every 
\m{c\in C} we have an exact sequence
\[0\lra\kn_c\lra\ko_{U_c^{(n)}}\ot\C^r\hfl{\phi_c}{}\ke_{c|U_c^{(n)}}\lra 0\]
(with \ \m{U_c=(X\times\{c\})\cap U}). Since the type of \m{\ke_c} is 
constant, \m{\phi_c} induces an isomorphism \ 
\m{\ko_{U_c}\ot\C^r\simeq\ke_{c|U_c}}. It follows that \m{\kn_c} is a quasi 
locally free sheaf on \m{U_c^{(n-1)}}, and its type is independent of $c$ (if 
\m{U_c} is non empty).

For every \m{c\in C} such that \m{U_c} is non empty, let \m{\ki_c} be the ideal 
sheaf of \m{X\times\{c\}\subset X_n\times C}. For every infinite subset 
\m{\Sigma\subset C}, we have \ 
\[\bigcap_{c\in\Sigma}\ki_c.(\ko_{U^{(n)}}\ot\C^r) \ = \ 0 \ . \] 
We have, if \m{U_c\not=\emptyset},
\[L^{n-1}.\kn \ \subset \ \ki_c.\kn \ \subset \ \ki_c.(\ko_{U^{(n)}}\ot\C^r)\]
(recall that \m{L^{n-1}} is the ideal sheaf of \m{X_{n-1}} in \m{X_n}). It 
follows that \ \m{L^{n-1}.\kn=0}, i.e. $\kn$ is a sheaf on \m{U^{(n-1)}}. 
Hence 
from the induction hypothesis, $\kn$ is quasi locally free. For every 
\m{(x,c)\in U},
\[{}^t\psi_{c,x}:\C^r\lra\kn_{c,x}^*\]
is surjective. It follows from corollary \ref{cor5} that \m{\ke_{|U^{(n)}}} 
is 
quasi locally free. This proves {\bf 1}, and {\bf 2}, {\bf 3} are immediate 
consequences of {\bf 1}.
\end{proof}

\sepprop

There is an easy converse of theorem \ref{theo4}:

\sepprop

\begin{subsub}\label{prop34}{\bf Proposition: } Let $\ke$ be a coherent sheaf 
on \m{X_n\times C}, \m{c\in C}, \m{x\in X} closed points, and $U$ a 
neighborhood of \m{(x,c)} in \m{X_n\times C} such that \m{\ke_{|U}} is quasi 
locally free. Then \m{\ke_{|U}} is flat on $C$.
\end{subsub}

\end{sub}

\sepsub

\Ssect{The Picard group of a primitive double scheme}{pic-sch}

(cf. \cite{dr11}, 7.)

Let \ \m{\boldsymbol{P}\subset\Pic(X)} \ be an irreducible component, such that 
some line bundle in $\boldsymbol{P}$ can be extended to \m{X_2}.
Let \m{\Pic^{\boldsymbol{P}}(X_2)} be the set of line bundles on \m{X_2} whose 
restriction to $X$ belongs to $\boldsymbol{P}$. If \ 
\m{\boldsymbol{P}=\Pic^0(X)} \ we will note \ 
\m{\Pic^{\boldsymbol{P}}(X_2)=\Pic^0(X_2)}. Let 
\m{\Gamma^0(X_2)\subset\Pic^0(X)} (resp. 
\m{\Gamma^{\boldsymbol{P}}(X_2)\subset\Pic^{\boldsymbol{P}}(X)}) be the set of 
line bundles that can be extended to \m{X_2}. Then \m{\Gamma^ 
{\boldsymbol{P}}(X_2)} is a smooth variety, isomorphic to 
\m{\Gamma^0(X_2)\subset\Pic^0(X)}, and \m{\Gamma^0(X_2)} is a subgroup of 
\m{\Pic^0(X)}.

The variety \m{\Pic^0(X_2)} is an algebraic group, the 
\m{\Pic^{\boldsymbol{P}}(X_2)} 
have a natural structure of smooth varieties, and the choice of an element 
\m{\Pic^{\boldsymbol{P}}(X_2)} defines in an obvious way an isomorphism \ 
\m{\Pic^0(X_2)\simeq\Pic^{\boldsymbol{P}}(X_2)}. The 
\m{\Pic^{\boldsymbol{P}}(X_2)} are the irreducible components of \m{\Pic(X_2)}.
Each \m{\Pic^{\boldsymbol{P}}(X_2)} has a natural structure of affine bundle on 
\m{\Gamma^{\boldsymbol{P}}(X_2)}, with associated vector bundle \ \m{\ko_X\ot 
H^1(X,L)}. In general this affine bundle is not banal (cf. \cite{dr11}, theorem 
9.2.1 when $X$ is a curve).

There is an open cover \m{(P_i)_{i\in I}} of \m{\Pic^{\boldsymbol{P}}(X_2)}, 
and for every \m{i\in I}, a {\em Poincaré bundle} \m{\L_i} on \ \m{X_2\times 
P_i} such that, with these data, \m{\Pic^{\boldsymbol{P}}(X_2)} is a {\em fine 
moduli space} in the sense of \ref{fin-mod}.

\end{sub}

\newpage

\section{Reduced Chern classes and reduced Hilbert polynomials}\label{RCCRHP}

\Ssect{Definitions}{RP_def}

If $Z$ is a projective scheme,  \m{\ko_Z(1)} an ample line bundle  and $E$ a 
coherent sheaf on $Z$, let \m{P_{\ko_Z(1)}(E)} be the Hilbert polynomial of $E$ 
with respect to \m{\ko_Z(1)}.

Let $X$ be a smooth, projective and irreducible variety. Let \m{X_n} be a
primitive multiple scheme of multiplicity $n$, with underlying smooth variety
$X$, and associated line bundle $L$ on $X$.

Let \m{\ko_X(1)} be an ample line bundle on $X$. For every coherent sheaf $E$
on $X$, let \ \m{c(E)\in A^*(X)} (or \m{H^*(X,\Z)}) be the total Chern class of
$E$.

Let $\ke$ be a coherent sheaf on \m{X_n}, and
\[\ke_n=0\subset\ke_{n-1}\subset\cdots\subset\ke_1\subset\ke_0=\ke\]
the first canonical filtration of $\ke$. Let
\[c_{red}(\ke) \ = \ \prod_{i=0}^{n-1}c(\ke_i/\ke_{i+1}) \ \in \ A^*(X) \
(\text{or}\ H^*(X,\Z)) \ , \qquad P_{red,\ko_X(1)}(\ke) = \
\sigg_{i=0}^{n-1}P_{\ko_X(1)}(\ke_i/\ke_{i+1}) \ . \]
We call \m{c_{red}(\ke)} the {\em total reduced Chern class of $\ke$}, and
\m{P_{red,\ko_X(1)}(\ke)} the {\em reduced Hilbert polynomial of $\ke$} (with 
respect to \m{\ko_X(1)}).

\sepprop

\begin{subsub}\label{rel_H} Relations with the usual Hilbert polynomial -- \rm
If \m{\ko_X(1)} can be extended to a line bundle \m{\ko_{X_n}(1)} on $Y$, then
\m{\ko_{X_n}(1)} is ample (prop. 4.6.1 of \cite{dr10}), and \
\m{P_{red,\ko_X(1)}(\ke)=P_{\ko_{X_n}(1)}(\ke)}. 
\end{subsub}

\sepprop

\begin{subsub}\label{prop1}{\bf Proposition:} {\bf 1 -- } Let \
\m{\kf_m=0\subset\kf_{m-1}\subset\cdots\subset\kf_1\subset\kf_0=\ke}
 \ be a
filtration such that, for \m{0\leq i<m}, \m{\kf_i/\kf_{i+1}} is 
supported on $X$. Then we have
\[c_{red}(\ke) \ = \ \prod_{i=0}^{m-1}c(\kf_i/\kf_{i+1}) \ , \qquad
P_{red,\ko_X(1)}(\ke) = \ \sigg_{i=0}^{m-1}P_{\ko_X(1)}(\kf_i/\kf_{i+1}) \ . 
\]
{\bf 2 -- } Let \
\m{0\to\ke'\to\ke\to\ke''\to 0} be an exact sequence of coherent sheaves on
\m{X_n}. Then we have
\[c(\ke) \ = c(\ke')c(\ke'') \quad \text{and} \quad P_{red,\ko_X(1)}(\ke) \ = \
P_{red,\ko_X(1)}(\ke')+P_{red,\ko_X(1)}(\ke'') \ . \]
\end{subsub}

\begin{proof} {\bf 1-} is an immediate consequence of proposition \ref{theo7}.
We now prove {\bf 2-}. According to proposition \ref{theo7} there exist
similar refinements \m{(\ke'_j)_{0\leq j\leq p}} of the filtration 
\m{\ke'\subset\ke} and \m{(\kf'_j)_{0\leq j\leq p}} of the first canonical 
filtration of $\ke$. There exists an integer \m{j_0} such that \m{0\leq 
j_0\leq 
p} and \m{\kf'_{j_0}=\ke'}. We have
\[c_{red}(\ke) \ = \ \prod_{j=0}^{p-1}c(\kf'_j/\kf'_{j+1}) \ , \]
and from {\bf 2-}
\[c_{red}(\ke') \ = \ \prod_{j=j_0}^{p-1}c(\kf'_j/\kf'_{j+1}) \ , \qquad
c_{red}(\ke'') \ = \ \prod_{j=0}^{j_0-1}c(\kf'_j/\kf'_{j+1}) \ . \]
Hence \ \m{c_{red}(\ke)=c_{red}(\ke')c_{red}(\ke)}, and similarly \
\m{P_{red,\ko_X(1)}(\ke)=P_{red,\ko_X(1)}(\ke')+P_{red,\ko_X(1)}(\ke'')}.
\end{proof}

\end{sub}

\sepsub

\Ssect{Invariance by deformation}{RP_surfa}

Let $C$ be an irreducible smooth curve and $\ke$ a coherent sheaf on 
\m{X_n\times C}, flat on $C$. The following result was pointed out by J\'anos 
Koll\'ar:

\sepprop

\begin{subsub}\label{theo8}{\bf Theorem: } Suppose that we have a filtration
\[0=\kg_m\subset\kg_{m-1}\subset\cdots\subset\kg_1\subset\kg_0=\ke\]
such that for \ \m{0\leq i<m}, \m{\kg_i/\kg_{i+1}} is supported on 
\m{X\times C}. Then there exists a filtration
\begin{equation}\label{equ24}
0=\kf_m\subset\kf_{m-1}\subset\cdots\subset\kf_1\subset\kf_0=\ke
\end{equation}
such that
\begin{enumerate}
\item[--] $\kf_1,\ldots,\kf_{m-1}$ and \ 
\m{\kf_0/\kf_1,\ldots,\kf_{m-2}/\kf_{m-1}} are flat on $C$. 
\item[--] $\kf_0/\kf_1,\ldots,\kf_{m-1}/\kf_{m}$ are supported on \m{X\times 
C}.
\item[--] There exists a finite subset \m{\Sigma\subset C} such that the 
two preceding filtrations restricted to \ \m{X_n\times (C\backslash\Sigma)} \ 
are the same.
\end{enumerate}
\end{subsub}
\begin{proof} Let \m{c_0\in C} and \m{\ki_{c_0}} its ideal sheaf. Let $z$ be a 
generator of the maximal ideal of \m{\ko_{C,c_0}}. Let \m{F=\kg_0/\kg_1}, and 
\m{G\subset F} be the subsheaf union of the subsheaves annihilated by the 
powers of \m{\ki_{c_0}}.

Then for every \m{x\in X}, \m{(F/G)_{(x,c_0)}} is a flat 
\m{\ko_{C,c_0}}-module: Let \m{\ov{u}\in (F/G)_{(x,c_0)}} and suppose that 
\m{z\ov{u}=0}. Let \m{u\in F_{(x,c_0)}} be over \m{\ov{u}}. Then \m{zu\in 
G_{(x,c_0)}}, so there exists an integer \m{k\geq 0} such that 
\m{z^k.zu=z^{k+1}u=0}. Hence \m{u\in G_{(x,c_0)}} and \m{\ov{u}=0}. 

Let \m{H\subset\ke} be the inverse image of $G$. From \cite{sga1}, 
Expos\'e IV, proposition 1.1, \m{H_{(x,c_0)}} is a flat 
\m{\ko_{C,c_0}}-module.

Since $X$ is projective, there is a neighborhood \m{C_0} of \m{c_0} such that 
\m{(F/G)_{|X\times C_0}} and \m{H_{|X_n\times C_0}} are flat on \m{C_0}. 
Note that $H$ coincide with \m{\kg_1} on \ \m{X_n\times(C_0\backslash\{c_0\})}. 
By taking othen \m{c_0} we can cover $C$ with a finite number of such open 
subsets \m{C_0} and finally 
obtain a subsheaf \m{\kg'_1\subset\ke} such that:
\begin{enumerate}
\item[--] $\kg_1\subset\kg'_1$ and $\kg'_1$ is flat on $C$.
\item[--] $\ke/\kg'_1$ is flat on $C$ and supported on \m{X\times C}.
\item[--] There is an open subset $\Sigma_1\subset C$ such that 
\m{\kg'_{1|X_n\times(C\backslash\Sigma_1)}=
\kg_{1|X_n\times(C\backslash\Sigma_1)}}.
\end{enumerate}
We can continue this process with \m{\kg'_1} instead of $\ke$, with the 
filtration \Nligne 
\m{0=\kg_m\subset\kg_{m-1}\subset\cdots\subset\kg_2\subset\kg'_1}, 
and so on. We finally obtain the filtration
\[0\subset\kf_m\subset\kf_{m-1}\subset\cdots\subset\kf_1\subset\kf_0=\ke\]
by subsheaves flat on $C$, as well as the quotients \m{\kf_i/\kf_{i+1}}, which 
are supported on \m{X\times C}, and this filtration coincide with the 
original one on \ \m{X_n\times (C\backslash\Sigma)}, where \m{\Sigma\subset C} 
is finite. But we have \m{\kf_m=0}, since it is flat on $C$ and 0 on 
\m{X_n\times(C\backslash\Sigma)}.
\end{proof}

\sepprop

\begin{subsub}\label{conj}{\bf Corollary: } The map 
\[\xymatrix@R=5pt{C\ar[r] & \Q[T]\\ c\fmaps[r] & P_{red,\ko_X(1)}(\ke_c)
}\]
is constant.
\end{subsub}
\begin{proof} We can view \m{X_n\times C} as a primitive multiple scheme, and \ 
\m{(X_n\times C)_{red}=X\times C}. Now, in theorem \ref{theo8}, for the first 
filtration we take the first canonical filtration of $\ke$. We obtain the 
filtration $(\ref{equ24})$. For every \m{c\in C}, we get a filtration of 
\m{\ke_c}
\[0=\kf_{m,c}\subset\kf_{m-1,c}\subset\cdots\subset\kf_{1,c}\subset\kf_{0,c}
=\ke_c \ , \]
such that for \m{0\leq i<m}, \m{\kf_{i,c}/\kf_{i-1,c}} is supported on 
\m{X\times C} and isomorphic to \m{(\kf_i/\kf_{i-1})_c}. Hence we have
\[P_{red,\ko_X(1)}(\ke_c) \ = \ 
\sigg_{i=0}^{m-1}P_{\ko_X(1)}\big((\kf_i/\kf_{i-1})_c\big) \ . \]
Since \m{\kf_i/\kf_{i-1}} is flat on $C$, \ \m{c\mapsto 
P_{\ko_X(1)}\big((\kf_i/\kf_{i-1})_c\big)} \ is constant, and so is \ 
\m{c\mapsto P_{red,\ko_X(1)}(\ke_c)}.
\end{proof}

\end{sub}

\sepsec

\section{Balanced sheaves on primitive multiple schemes}\label{CSPMS}

\Ssect{Balanced sheaves}{ext_hi}

Let $X$ be a smooth and irreducible variety. Let \m{X_n} be a primitive 
multiple scheme of multiplicity $n$, with underlying smooth variety $X$, and 
associated line bundle $L$ on $X$.Let $\ke$ be a coherent sheaf on \m{X_n}. 
From \ref{prop_filt} we have canonical morphisms
\[\lambda_{i-1}\ot I_{L^{i-2}}:G^{(i)}(\ke)\ot L^{i-1}\lra G^{(i-1)}(\ke)\ot 
L^{i-2} \qquad \text{for} \ \ 2\leq i\leq n \ , \]
which is injective, and
\[\mu_i\ot I_{L^{n-2-i}}:G_i(\ke)\ot L^{n-1-i}\lra G_{i+1}(\ke)\ot L^{n-2-i} 
\qquad \text{for} \ \ 0\leq i\leq n-2 \ , \]
which is surjective. Let 
\[\boldsymbol{\lambda}(\ke) \ = \lambda_1\circ(\lambda_2\ot I_L)\circ\cdots\circ
(\lambda_{n-1}\ot I_{L^{n-2}}):G^{(n)}(\ke)\ot L^{n-1}\lra G^{(1)}(\ke) \ , \]
\[\boldsymbol{\mu}(\ke) \ = \ \mu_{n-2}\circ(\mu_{n-3}\ot I_L)\circ\cdots\circ
(\mu_0\ot I_{L^{n-2}}):G_0(\ke)\ot L^{n-1}\lra G_{n-1}(\ke) \ . \]

\sepprop

We will use the easy following lemma:

\sepprop

\begin{subsub}\label{lem11}{\bf Lemma: } Let $A$ be a commutative ring, 
\m{k\geq 3} an integer, \m{M_1,\ldots,M_k}\Nligne $A$-modules, and \ 
\m{f_i:M_i\to M_{i+1}}, \m{1\leq i<k}, injective (resp. surjective) morphisms. 
Then \ \m{f_{k-1}\circ\cdots f_1:M_1\to M_k} \ is surjective (resp. injective) 
if and only \m{f_1,\ldots,f_k} are surjective (resp. injective).
\end{subsub}

\sepprop

\begin{subsub}\label{def1}{\bf Definition:} We say that $\ke$ is {\em balanced} 
if \m{\boldsymbol{\lambda}(\ke)} is surjective.
\end{subsub}

\sepprop

From lemma \ref{lem11} and proposition \ref{prop15}, we have

\sepprop

\begin{subsub}\label{prop16}{\bf Proposition: } The following properties are 
equivalent
\begin{enumerate}
\item[(i)] $\ke$ is balanced.
\item[(ii)] $\lambda_1,\ldots,\lambda_{n-1}$ are surjective
\item[(iii)] $\Gamma_1(\ke)=\cdots=\Gamma_{n-1}(\ke)=0$ .
\item[(iv)] $\Gamma^{(1)}(\ke)=\cdots=\Gamma^{(n-1)}(\ke)=0$ .
\item[(v)] $\mu_0,\ldots,\mu_{n-2}$ are injective.
\item[(vi)] $\boldsymbol{\mu}(\ke)$ is injective.
\end{enumerate}
If $\ke$ is balanced, then \m{\lambda_1,\ldots,\lambda_{n-1}},
\m{\mu_0,\ldots,\mu_{n-2}}, $\boldsymbol{\lambda}$ and $\boldsymbol{\mu}$ are 
isomorphisms.
\end{subsub}

\sepprop

\begin{subsub}\label{prop17}{\bf Proposition: } Let $\ke$ be a coherent sheaf 
on \m{X_n}. Then $\ke$ is balanced if and only if \ \m{\ke_i=\ke^{(n-i)}} \ for 
\m{1\leq i\leq n}.
\end{subsub}
\begin{proof} If \ \m{\ke_i=\ke^{(n-i)}} \ for \m{1\leq i\leq n}, we have \ 
\m{\lambda_i=\mu_{n-i-1}}, hence \m{\lambda_i} and \m{\mu_{n-i-1}} are 
isomorphisms, and $\ke$ is balanced.

Conversely, suppose that $\ke$ is balanced.
For \m{1\leq j\leq n}, let \ \m{\beta_j:\ke_{n-j}/\ke_{n-j+1}\to
\ke^{(j)}/\ke^{(j-1)}} \ be the morphism induced by the inclusion \ 
\m{\ke_{n-j}\subset\ke^{(j)}}.

We have a commutative diagram
\xmat{\ke/\ke_1\ot L^{n-1}\ar[r]^-{\boldsymbol{\mu}}\flon[d]^\alpha &
\ke_{n-1}\flinc[d]^{\beta_1}\\
\ke/\ke^{(n-1)}\ar[r]^-{\boldsymbol{\lambda}} & \ke^{(1)}
}
where the surjective morphism $\alpha$ is induced by the inclusion \ 
\m{\ke_1\subset\ke^{(n-1)}} \ and $\beta$ is the inclusion. Since 
$\boldsymbol{\lambda}$ and $\boldsymbol{\mu}$ are isomorphisms, so are $\alpha$ 
and \m{\beta_1}. For \m{2\leq i\leq n-2} we have a commutative diagram
\xmat{(\ke_{n-i}/\ke_{n-i+1})\ot L\ar[rr]^-{\mu_{n-i}}\ar[d]^{\beta_i} & &
\ke_{n-i+1}/\ke_{n-i+2}\ar[d]^{\beta_{i-1}}\\
(\ke^{(i)}/\ke^{(i-1)})\ot L\ar[rr]^-{\lambda_{i-1}} & & \ke^{(i-1)}/\ke^{(i-2)}
}
Using the fact that \m{\lambda_{i-1}}, \m{\mu_{n-i}}, \m{\beta_1} are 
isomorphisms we see by induction on $i$ that \m{\beta_i} is an isomorphism. 
Again by induction on $i$ it is then easy to see that \ \m{\ke_i=\ke^{(n-i)}}.
\end{proof} 

\sepprop

\begin{subsub}\label{ex2}{\bf Examples: }\rm -- Vector bundles on \m{X_n} are 
balanced sheaves.

-- Recall the a sheaf of ideals $\ki$ on \m{X_n} is called {\em regular} if for 
every closed point \m{x\in X_n}, \m{\ki_x} is generated by the elements of a 
regular sequence in \m{\ko_{X_n,x}}. It follows from proposition \ref{prop21} 
that a regular ideal is a balanced sheaf.
\end{subsub}

\sepprop

\begin{subsub}\label{bal2} The case of double schemes -- \rm
We suppose that \m{n=2}. Let $\ke$ be a coherent sheaf of 
\m{X_2}. Then we have \ \m{\ke_1\subset\ke^{(1)}}. Let \ \m{F\subset\ke} \ be a 
subsheaf and \ \m{E=\ke/F}. Then $F$ is supported on $X$ if and only if \ 
\m{F\subset\ke^{(1)}}, and in this case $E$ is supported on $X$ if and only 
if \ \m{\ke_1\subset F}.
\end{subsub}

Suppose that \ \m{\ke_1\subset F\subset\ke^{(1)}}. The canonical morphism \ 
\m{\ke\ot L\to\ke} \ induces \ \m{\phi_F:E\ot L\to F}, and \m{\phi_F} is 
injective (resp. surjective) if and only if \m{F=\ke^{(1)}} (resp. \m{F=\ke_1}).
The sheaf $\ke$ is supported on $X$ if and only if \m{\phi_F=0}.

If $\ke$ is balanced then  \ \m{\ke\ot L\to\ke} \ induces a canonical 
isomorphism \ \m{\ke_1\simeq\ke_{|X}\ot L}.

Let \m{F\subset\ke} be a subsheaf supported on $X$. If \m{\phi_F} is an 
isomorphism, then $\ke$ is balanced, and \m{F=\ke_1}. It follows that

\sepprop

\begin{subsub}\label{lem8}{\bf Lemma: } Let $E$, $\ke$ be coherent sheaves on 
$X$, \m{X_2} respectively. Suppose that $E$ is simple, and that we have an 
exact sequence \ \m{0\to E\ot L\to\ke\to E\to 0}. If $\ke$ is not supported 
on $X$, then $\ke$ is a balanced sheaf and \m{\ke_{|X}=E}.
\end{subsub}

\end{sub}

\sepsub

\Ssect{Flat families of balanced sheaves on primitive double schemes}{fam_bal_1}

We suppose that \m{n=2} and that $X$ is projective. Let $C$ be an irreducible 
smooth curve. We can view \m{X_2\times C} as a primitive double scheme.

Let $Z$ be a scheme over $\C$. We can extend the definitions of 
\ref{QF} and of balanced sheaves to \m{X_2\times Z}.

Let $\kf$ be a coherent sheaf on \m{X_2\times Z}, flat on $Z$. We can consider 
the property of being balanced for $\ke$, or for the fibers \m{\kf_z} (where 
$z$ is a closed point of $Z$). We will study the relations between the two. In 
particular we obtain that the property of being balanced is an open property.

\sepprop

\begin{subsub}\label{theo5}{\bf Theorem: }
{\bf 1 -- } Let $\ke$ be a coherent sheaf on \m{X_2\times Z}, flat on $Z$. 
Suppose that for some closed point \m{z\in Z}, \m{\ke_z} is balanced. Then 
there exists a neighborhood $V$ of $z$ such that \m{\ke_{|X_2\times V}} is 
balanced, \m{\ke_{1|X_2\times V}}, \m{\ke_{|X\times V}} are flat on $V$ and \ 
\m{\ke_{1,v}=\ke_{v,1}}, \m{(\ke_{|X\times V})_v=\ke_{v|X}} \ for every closed 
point \m{v\in V}.

{\bf 2 -- } Let $\ke$ be a coherent sheaf on \m{X_2\times C}, flat on $C$. 
Suppose that $\ke$ is balanced. Then for every \m{c\in C}, \m{\ke_c} is 
balanced.
\end{subsub}
\begin{proof} We first prove {\bf 1}. Let \m{x_0\in X}, and $U$ a neighborhood 
of \m{x_0} such \m{L_{|U}} can be extended to a line bundle $\L$ on 
\m{U^{(2)}}. Let \ \m{p:U^{(2)}\times Z\to U^{(2)}} \ be the projection. We 
have canonical exact sequences
\[0\lra(\ke^{(1)})_{|U^{(2)}\times Z}\ot p^*(\L)\lra\ke_{|U^{(2)}\times Z}\ot 
p^*(\L)\lra\ke_{1|U^{(2)}\times Z}\lra 0 \ , \]
\[0\lra\ke_1\lra\ke\lra\ke_{|X\times Z}\lra 0 \ . \]
From the first one we deduce the exact sequence on $U$
\[0\lra\Tor^1_{\ko_{X_2\times Z}}(\ke_1,\ko_{X_2\times\{z\}})\lra
(\ke^{(1)})_z\ot\L\lra\ke_z\ot\L\lra\ke_{1,z}\lra 0 \ . \]
It follows from proposition \ref{prop30} that \ \m{\ke_{1,z}=\ke_{z,1}} \ on 
$U$. In particular the morphism \m{\ke_{1,z}\to\ke_z} induced by the 
inclusion 
\m{\ke_1\subset\ke} is injective on \m{U^{(2)}}. From the second 
canonical exact sequence and\cite{sga1}, Expos\'e IV, corollaire 5.7, 
\m{(\ke_{|X\times Z})_{(x,z)}} is a flat \m{\ko_{Z,z}}-module, for every \ 
\m{x\in U}. By \cite{sga1}, Expos\'e IV, Theorem 6.10, there is an 
open neighborhood \m{W_0} of $z$ such that \m{\ke_{|U\times W_0}} is 
flat on \m{W_0}, and from \cite{sga1}, Expos\'e IV, proposition 1.1, 
\m{\ke_{1|U\times W_0}} is flat on \m{W_0}. We can then take a finite number of 
points \m{x_0\in X} such that the union of corresponding open subsets 
\m{U\subset X} is X. Taking the intersection $W$ of the corresponding open 
subsets \m{W_0\subset Z}, we see that \m{\ke_{1|X\times W}} \m{\ke_{|X\times 
W}} are flat on $W$.

From proposition \ref{prop30} we have \ \m{(\ke_{|X\times Z})_z=\ke_{z|X}}. 
Hence the canonical surjective morphism \ \m{\Phi:\ke_{|X\times Z}\to\ke_1\ot 
L^*} \ is an isomorphism on \ \m{X\times\{z\}}. Let \ 
\m{\kn=\ker(\Phi)}. Then \m{\kn_{|X\times W}} is flat on $W$. Let 
\m{\mm_z\subset\ko_{Z,z}} be the maximal ideal, and \m{x\in X}. Then since 
\m{\kn_z=0}, we have \m{\kn_{(x,z)}=\mm_z\kn_{(x,z)}}, and \m{\kn_{(x,z)}=0}. 
To prove {\bf 1} we take \m{V\subset W} such that \m{X\times V} does not meet 
the support of $\kn$.

Now we prove {\bf 2}. Let \m{c\in C}, \m{x\in X}, $t$ a local section of 
$L$, defined at $x$ and that generates $L$ around $x$, and $z$ a generator of 
\m{\mm_c}. We must prove that \ \m{(\ke_c)_{1,x}=(\ke_c)^{(1)}_x}. Let 
\m{e_c\in(\ke_c)^{(1)}_x}, i.e. \m{te_c=0}. Let \m{e\in\ke_{(x,c)}} over 
\m{e_c}. Then there exists \m{u\in\ke_{(x,c)}} such that \ \m{te=zu}. We have \ 
\m{z.tu=t.zu=t^2e=0}, hence, since \m{\ke_{(x,c)}} is a flat 
\m{\ko_{C,c}}-module, \m{tu=0}. Since $\ke$ is balanced, we can write \ 
\m{u=tv}, with \m{v\in\ke_{(x,c)}}. We have \ \m{t(e-zv)=0}, and since $\ke$ is 
balanced, there exists \m{f\in\ke_{(x,c)}} such that \ \m{e-zv=tf}. If \m{f_c} 
is the image of $f$ in \m{\ke_{c,x}}, we have \ \m{e_c=tf_c}, and \ 
\m{e_c\in(\ke_c)_{1,x}}. Hence \ \m{(\ke_c)^{(1)}_x\subset(\ke_c)_{1,x}}. This 
proves {\bf 2}.
\end{proof}

\sepprop

\begin{subsub}\label{cor6}{\bf Corollary: } Let $\ke$ be a coherent sheaf on 
\m{X_2\times C}, flat on $C$. Suppose that for some \m{c\in C}, \m{\ke_c} is 
balanced. Then there exists a neighborhood $V$ of $c$ such that for every 
\m{v\in V}, \m{\ke_v} is balanced, \m{\ke_{1|X_2\times V}}, \m{\ke_{|X\times 
V}} are flat on $V$ and \ \m{\ke_{1,v}=\ke_{v,1}}, \m{(\ke_{|X\times 
V})_v=\ke_{v|X}} \ for every \m{v\in V}.
\end{subsub}

\end{sub}

\sepsub

\Ssect{Flat families of balanced sheaves on primitive double surfaces}{fam_bal}

We suppose that \m{n=2}, $X$ is projective, and \ \m{\dim(X)=2}. Let $S$ be a 
smooth variety and \m{{\bf F}} be a family of torsion free sheaves on $X$, 
parameterized by $S$ and flat on $S$. We will prove in corollary \ref{coro4} 
that, with suitable conditions, the set
\[\{s\in S; {\bf F}_s\ \text{can be extended to a balanced sheaf on} \ X_2\}\] 
is closed is $S$. In \cite{dr11}, moduli spaces of vector bundles on \m{X_2} 
are constructed from moduli spaces of vector bundles on $X$ (by trying to 
extend the vector bundles on $X$ to \m{X_2}). Corollary \ref{coro4} should help 
to extend these constructions to moduli spaces of torsion free sheaves on $X$, 
and obtain moduli spaces of balanced sheaves on \m{X_2}.

Let \m{{\bf E}} be another family of torsion free sheaves on $X$, parameterized 
by $S$ and flat on $S$. Let 
\[f_0:\F_0\lra{\bf F}\]
be a surjective morphism, where \m{\F_0} is a vector bundle on \m{X\times S} 
such that for every \m{s\in S} and positive integer $i$ we have \ 
\m{h^i(X,\F_{0s}^*\ot{\bf E}_s)=0}. Let \ \m{X_0=\ker(f_0)}. By \ref{loc_free}, 
\m{X_0} is a flat family of torsion free sheaves, and by \cite{ha3}, prop. 1.1, 
for every \m{s\in S}, \m{X_{0,s}} is a reflexive sheaf, and since 
\m{\dim(X)=2}, \m{X_{0,s}} is locally free. Hence \m{X_0} is locally free.

Let \m{s\in S}. From the exact sequence \ \m{0\to X_{0,s}\to\F_{0,s}\to{\bf F}_s
\to 0}, we have an exact sequences
\[0\to\Hom({\bf F}_s,{\bf E}_s)\lra\Hom(\F_{0,s},{\bf 
E}_s)\lra\Hom(X_{0,s},{\bf E}_s)\lra\qquad\qquad\qquad\]
\[\qquad\qquad\qquad \Ext^1_{\ko_X}({\bf F}_s,{\bf 
E}_s)\lra\Ext^1_{\ko_X}(\F_{0,s},{\bf E}_s)=\nsp , \]
and
\[0\to\Hom({\bf F}_s,{\bf E}_s)\lra\Hom(\F_{0,s},{\bf 
E}_s)\lra\Hom(X_{0,s},{\bf E}_s)\lra\qquad\qquad\qquad\]
\[\qquad\qquad\qquad\Ext^1_{\ko_{X_2}}({\bf F}_s,{\bf E}_s)\hfl{\delta_s}{}
\Ext^1_{\ko_{X_2}}(\F_{0,s},{\bf E}_s)\hfl{\theta}{}
\Ext^1_{\ko_{X_2}}(X_{0,s},{\bf E}_s) \ . \]
Hence we have an exact sequence
\[0\lra\Ext^1_{\ko_X}({\bf F}_s,{\bf E}_s)\lra\Ext^1_{\ko_{X_2}}({\bf F}_s,{\bf 
E}_s)\hfl{\delta_s}{}\Ext^1_{\ko_{X_2}}(\F_{0,s},{\bf E}_s)\hfl{\theta}{}
\Ext^1_{\ko_{X_2}}(X_{0,s},{\bf E}_s) \ . \]

\sepprop

\begin{subsub}\label{lem10}{\bf Lemma: } Let $A$ (resp. B) be a vector bundle 
(resp. a coherent sheaf) on $X$. Then there is a canonical functorial 
isomorphism
\[\EExt^1_{\ko_{X_2}}(A,B) \ \simeq \ \HHom(A\ot L,B) \ . \]
\end{subsub}
The proof is similar to that of lemma 4.6.2 of \cite{dr11}.

\sepprop

By the Ext spectral sequence and the fact that \ \m{H^i(X,\F_{0,s}^*\ot{\bf 
E}_s)=\nsp} \ for \m{i=1,2}, the canonical map
\[\Ext^1_{\ko_{X_2}}(\F_{0,s},{\bf E}_s)\lra 
H^0(\EExt^1_{\ko_{X_2}}(\F_{0,s},{\bf E}_s))=\Hom(\F_{0,s}\ot L,{\bf E}_s)\]
is an isomorphism.

Hence we have a commutative diagram
\xmat{\Ext^1_{\ko_{X_2}}(\F_{0,s},{\bf E}_s)\ar[r]^-\theta\fleq[d] & 
\Ext^1_{\ko_{X_2}}(X_{0,s},{\bf E}_s)\ar[d]\\
\Hom(\F_{0,s}\ot L,{\bf E}_s)\ar[r]^-\alpha & \Hom(X_{0,s}\ot L,{\bf E}_s)
}
(where $\alpha$ is induced by \ \m{X_{0s}\subset\F_{0s}}).
We have \ \m{\ker(\alpha)\simeq\Hom({\bf F}_s\ot L,{\bf E}_s)}, hence \Nligne 
\m{\imm(\delta_s)=\ker(\theta)\subset\Hom({\bf F}_s\ot L,{\bf E}_s)}, and we 
have an exact sequence
\[0\lra\Ext^1_{\ko_X}({\bf F}_s,{\bf E}_s)\lra\Ext^1_{\ko_{X_2}}({\bf F}_s,{\bf 
E}_s)\hfl{\delta_s}{}\Hom({\bf F}_s\ot L,{\bf E}_s) \ . \]
We suppose now that
\begin{enumerate}
\item[--] ${\bf E}={\bf F}\ot p_X^*(L)$, where $p_X:X\times S\to X$ is the 
projection.
\item[--] For every $s\in S$, ${\bf F}_s$ is simple.
\item[--] $\dim(\Ext^1_{\ko_X}({\bf F}_s,{\bf F}_s\ot L))$ is independent of 
$s\in S$.
\end{enumerate}
It follows that we have an exact sequence
\[0\lra\Ext^1_{\ko_X}({\bf F}_s,{\bf F}_s\ot L)\lra\Ext^1_{\ko_{X_2}}({\bf 
F}_s,{\bf F}_s\ot L)\hfl{\delta_s}{}\C \ . \]

\sepprop

\begin{subsub}\label{prop14}{\bf Proposition: } Let \m{s\in S}, \m{\sigma\in 
\Ext^1_{\ko_{X_2}}({\bf F}_s,{\bf F}_s\ot L)}, and \ \m{0\to{\bf F}_s\ot L\to\ke
\to{\bf F}_s\to 0} \ the corresponding extension. Then $\ke$ is a balanced 
sheaf if and only if \ \m{\delta_s(\sigma)\not=0}.
\end{subsub}
\begin{proof} If  \m{\delta_s(\sigma)=0} then \ \m{\sigma\in\Ext^1_{\ko_X}({\bf 
F}_s,{\bf F}_s\ot L)}, hence $\ke$ is supported on $X$ and cannot be 
balanced. Conversely assume that \m{\delta_s(\sigma)\not=0}. By lemma 
\ref{lem8} it suffices to prove that $\ke$ is not supported on $X$, which is 
obvious since otherwise we would have \m{\sigma\in\Ext^1_{\ko_X}({\bf 
F}_s,{\bf F}_s\ot L)} and \m{\delta_s(\sigma)=0}.
\end{proof}

\sepprop

Let \ \m{n=\dim(\Ext^1_{\ko_X}({\bf F}_s,{\bf F}_s\ot L))}. It follows that 
\m{{\bf F}_s} can be extended to a balanced sheaf on \m{X_2} if and only if \ 
\m{\dim(\Ext^1_{\ko_{X_2}}({\bf F}_s,{\bf F}_s\ot L))=n+1}.

\sepprop

\begin{subsub}\label{coro4}{\bf Corollary: } The set \m{\{s\in S; {\bf F}_s\ 
\text{can be extended to a balanced sheaf}\}} is closed is $S$.
\end{subsub}
\begin{proof} This follows from the upper semicontinuity of the map \ 
\m{s\mapsto\dim(\Ext^1_{\ko_{X_2}}({\bf F}_s,{\bf F}_s\ot L))} .
\end{proof}

\end{sub}

\sepsec

\section{Extensions of ideal sheaves on primitive double surfaces}\label{ext_id}

Let $X$ be a complex smooth projective surface, $L$ a line bundle on $X$ and 
\m{X_2} a primitive double scheme, with underlying smooth variety $X$ and 
associated line bundle $L$. Let \m{P\in X} be a closed point.

\sepprop

{\em Notations -- } Let \m{\ko_1=\ko_{X,P}}, \m{\ko_2=\ko_{X_2,P}},
\m{\ki=\ki_{X,X_2,P}}, and $\mm$ the maximal ideal of \m{\ko_1}. 
We have \ \m{\ki\simeq\ko_1} \ and \ \m{\mm\ot_{\ko_1}\ki\simeq\mm\ki} .
Let $\ov{t}$ be a generator of $\ki$.

-- If $R$ is a ring and $n$ is a positive integer, \m{nR} will denote the 
direct sum of $n$ copies of $R$.

-- elements of \m{\ko_1} will be denoted by 
greek letters: \m{\alpha,\beta,\gamma,\ldots}

-- elements of \m{\ko_2} will be denoted by latin letters: 
\m{a,b,c,\ldots}

-- for the restrictions to $X$ we will use the index $_0$: 
\m{a_0,b_0,c_0,\ldots}

-- elements of the ideal \m{\ki} will be denoted by overscored 
letters: \m{\ov{a},\ov{\alpha},\ldots}

\sepsub

\Ssect{Ideals of \m{\ko_2}}{O2_id}

\begin{subsub}{\bf Definition: }\label{def3} \rm
Let $\J_P$ the set of ideals $\BS{J}$ of \m{\ko_2} such that there exist 
generators $x$, $y$ of $\BS{J}$ such that \ \m{\mm=(x_0,y_0)}.

Let \m{x,y\in\ko_2} be such that \m{x_0,y_0} are generators of $\mm$. Then we 
have 
\[\J_P \ = \ \{(x+\ov{A},y+\ov{B}) \ ; \ \ov{A},\ov{B}\in\ki \} \ . \]
\end{subsub}

\sepprop

We will first see how the elements \m{(x+\ov{A},y+\ov{B})} of $\J_P$ depend on 
the parameters \m{\ov{A},\ov{B}}.

\sepprop

\begin{subsub}\label{prop23}{\bf Proposition: } Let 
\m{\ov{A},\ov{B},\ov{A'},\ov{B'}\in\ki}, and $\tau$ (resp. \m{\tau'}) 
the image of \ \m{-y_0\ov{A}+x_0\ov{B}} (resp. \ \m{-y_0\ov{A'}+x_0\ov{B'}}) in 
\ \m{(\mm/\mm^2)\ot_{\ko_1}\ki}. The following properties are equivalent:
\begin{enumerate}
\item[(i)]$(x+\ov{A},y+\ov{B})=(x+\ov{A'},y+\ov{B'})$ .
\item[(ii)] $\tau=\tau'$ .
\item[(iii)] $\ov{A'}-\ov{A}\in\mm\ki$ \ and \ $\ov{B'}-\ov{B}\in\mm\ki$ .
\item[(iv)] the $\ko_2$-modules $(x+\ov{A},y+\ov{B})$ and 
$(x+\ov{A'},y+\ov{B'})$ are isomorphic.
\end{enumerate}
\end{subsub}

\begin{proof} Suppose that (iii) is true. Then 
\begin{equation}\label{equ7}
y_0(\ov{A'}-\ov{A}) \ = \ x_0(\ov{B'}-\ov{B}) \ = \ 0 \quad (\text{mod.} \
\mm^2\ki) \ , \end{equation}
hence \m{\tau=\tau'} and (ii) is true. Conversely, if (ii) is true, we have
\Nligne \m{y_0(\ov{A'}-\ov{A})=x_0(\ov{B'}-\ov{B}) \ (\text{mod.} \ \mm^2\ki)}.
Hence we can write
\[y_0(\ov{A'}-\ov{A}) \ = \ x_0(\ov{B'}-\ov{B})+\ov{t}\phi \ , \]
with \m{\phi\in\mm^2}. Suppose that \ \m{\ov{A'}-\ov{A}=a+\alpha}, 
\m{\ov{B'}-\ov{B}=b+\beta}, with \m{a,b\in\C}, \m{\alpha,\beta\in\mm}. We have
\[\big(y_0a-x_0b+[y_0\alpha-x_0\beta-\phi]\big)\ov{t} \ = \ 0 \ , \]
hence \ \m{y_0a-x_0b+[y_0\alpha-x_0\beta-\phi]=0}. Since \ 
\m{y_0\alpha-x_0\beta-\phi\in\mm^2}, we have \ \m{a=b=0}, and (iii) is true.

Now we prove that (i) implies (ii). Suppose that (i) is true.
There exist \m{\lambda,\mu,\epsilon,\rho\in\ko_2} such that
\begin{equation}\label{equ6}
x+\ov{A'} \ = \ \lambda(x+\ov{A})+\mu(y+\ov{B}) \ , \quad
y+\ov{B'} \ = \ \epsilon(x+\ov{A})+\rho(y+\ov{B}) \ .
\end{equation}
We have \ \m{\lambda_0x_0+\mu_0y_0=x_0}, hence there is 
some \m{\Psi\in\ko_2} such that \ \m{\lambda_0=1-\Psi y_0} \ and \ 
\m{\mu_0=\Psi x_0}. Similarly there exists \m{\theta\in\ko_2} such that \ 
\m{\epsilon_0=\theta y_0} \ and \ \m{\rho_0=1-\theta x_0}. Hence we can write
\[\lambda=1-\psi y+\ov{w} \ , \quad \mu=\psi x+\ov{z} \ , \quad
\epsilon=\theta y+\ov{u} \ , \quad \rho=1-\theta x+\ov{v} \ , \]
with \ \m{\ov{w},\ov{z},\ov{u},\ov{v}\in\ki}. It follows that, with \ 
\m{\ov{\tau}=-y_0\ov{A}+x_0\ov{B}}, \m{\ov{\tau'}=-y_0\ov{A'}+x_0\ov{B'}},
\[\ov{\tau'}-\ov{\tau} \ = \ -x_0y_0\ov{w}-y_0^2\ov{z}-y_0\psi_0\ov{\tau}
+x_0^2\ov{u}+x_0y_0\ov{v}-x_0\theta_0\ov{\tau} \ , \]
hence \m{\tau=\tau'} and (ii) is true.

Now we prove that (iii) implies (i). Suppose that \ \m{\ov{A'}-\ov{A}\in\mm\ki} 
\ and \ \m{\ov{B'}-\ov{B}\in\mm\ki}. So we can write \ 
\m{\ov{A'}=\ov{A}+x_0\ov{u}+y_0\ov{v}}, with \ \m{\ov{u},\ov{v}\in\ki}. Hence
\[x+\ov{A'} \ = \ x+\ov{A}+x_0\ov{u}+y_0\ov{v} \ = \ (1+\ov{u})(x+\ov{A})+
(y+\ov{B})\ov{v} \ , \]
(since \m{\ov{A}\ov{u}=\ov{B}\ov{v}=0}) so \ 
\m{x+\ov{A'}\in(x+\ov{A},y+\ov{B})}. Similarly 
\ \m{y+\ov{B'}\in(x+\ov{A},y+\ov{B})}. Hence \Nligne 
\m{(x+\ov{A'},y+\ov{B'})\subset(x+\ov{A},y+\ov{B})}. In the same way \ 
\m{(x+\ov{A},y+\ov{B})\subset(x+\ov{A'},y+\ov{B'})} \ and (i) is true.

Suppose now that (iii) is true. Then we have proved that (i) is true, and so is 
(iv). Conversely, suppose that (iv) is true. Let \ 
\m{\sigma:(x+\ov{A},y+\ov{B})\to(x+\ov{A'},y+\ov{B'})} \ be an isomorphism. It 
induces an automorphism of the \m{\ko_1}-module \m{(x_0,y_0)}, which is the 
multiplication by some invertible element of \m{\ko_1}. Hence with respect to 
the generators \m{x+\ov{A}}, \m{y+\ov{B}} of \m{(x+\ov{A},y+\ov{B})}, and
\m{x+\ov{A'}}, \m{y+\ov{B'}} of \m{(x+\ov{A'},y+\ov{B'})}, $\sigma$ is 
represented by a matrix \m{\begin{pmatrix}\alpha+\ov{u} & \ov{v}\\ \ov{w} & 
\alpha+\ov{z}\end{pmatrix}}, with \m{\alpha\in\ko_2} invertible, 
\m{\ov{u},\ov{v},\ov{w},\ov{z}\in\ki}. Replacing $\sigma$ with 
\m{\alpha^{-1}\sigma}, we can assume that \m{\alpha=1}. The equation
\[(y+\ov{B'})\sigma(x+\ov{A'}) \ = \ (x+\ov{A'})\sigma(y+\ov{B'})\]
gives
\[x_0[y_0(\ov{u}-\ov{z})+\ov{B'}-\ov{B}+x_0\ov{w}] \ = \
y_0[y_0\ov{v}+\ov{A'}-\ov{A}]\]
in $\ki$. It follows that \ \m{ \ov{A'}-\ov{A}\in\mm\ki}, and similarly \ 
\m{\ov{B'}-\ov{B}\in\mm\ki}. So (iii) is true.
\end{proof}

\end{sub}

\sepsub

\Ssect{Extensions of ideals of \m{\ko_{X,P}}}{ext_id_m}

\begin{subsub}\label{theo2}{\bf Theorem: } The balanced \m{\ko_2}-modules 
$M$ such that \ \m{M\ot_{\ko_2}\ko_1\simeq\mm} \ are the ideals of $\J_P$.
\end{subsub}
\begin{proof}
The ideals of $\J_P$ are balanced modules by proposition \ref{prop21}, and if 
\m{\kj\in\J_P}, then \ \m{\kj\ot_{\ko_2}\ko_1\simeq\mm}. Conversely we 
will study the extensions of \m{\ko_2}-modules
 \[0\lra\mm\ot_{\ko_1}\ki\to M\lra\mm\lra 0 \ , \]
and prove that when $M$ is balanced, then it is an ideal of $\J_P$.

We have the following free resolution of $\mm$ as a \m{\ko_2}-module:
\begin{equation}\label{equ3}\xymatrix{3\ko_2\ar[r]^-{\phi_2} & 
3\ko_2\ar[r]^-{\phi_1} & 2\ko_2\ar[r]^-{\phi_0} & \mm\ar[r] & 0 \ ,}
\end{equation}
where \m{\phi_0}, \m{\phi_1}, \m{\phi_2} are represented by the matrices
\[\begin{pmatrix}x & y\end{pmatrix} \ , \quad \begin{pmatrix}y & \ov{t} & 0\\
-x & 0 & \ov{t}\end{pmatrix} \ , \quad \begin{pmatrix}\ov{t} & 0 & 0\\
-y & \ov{t} & 0 \\ x & 0 & \ov{t}\end{pmatrix} \]
respectively. We have
\[K \ = \ker(\phi_0) \ = \{(ey+\ov{\gamma},-ex+\ov{\delta}); \ e\in\ko_2, \ 
\ov{\gamma},\ov{\delta}\in\ki\} \ . \]

\sepsubsub

{\bf Step 1. }{\em Morphisms from $K$ to \ \m{\mm\ot_{\ko_1}\ki=\mm\ki} -- }
Recall 
that $\ki$ is isomorphic to \m{\ko_1}, so \ \m{\mm\ot_{\ko_1}\ki\simeq\mm}. 
These 
morphisms are the \m{\Theta_{\ov{\tau},\rho}}, 
\m{\ov{\tau}\in\mm\ot_{\ko_1}\ki}, 
\m{\rho\in\ko_1}, with
\[\Theta_{\ov{\tau},\rho}(ey+\ov{\gamma},-ex+\ov{\delta}) \ = \
e_0\ov{\tau}+\rho(x_0\ov{\gamma}+y_0\ov{\delta}) \ . \]

\sepsubsub

{\bf Step 2. }{\em Description of \ \m{\Ext^1_{\ko_2}(\mm,\mm\ot_{\ko_1}\ki)} 
-- } From $(\ref{equ3})$ follows the complex
\xmat{2(\mm\ot_{\ko_1}\ki)\fleq[d]\ar[r]^-{\Psi_1} & 
3(\mm\ot_{\ko_1}\ki)\fleq[d]
\ar[r]^-{\Psi_2} & 3(\mm\ot_{\ko_1}\ki)\fleq[d]\\
\Hom(2\ko_2,\mm\ot_{\ko_1}\ki) & \Hom(3\ko_2,\mm\ot_{\ko_1}\ki) & 
\Hom(3\ko_2,\mm\ot_{\ko_1}\ki) 
\ , }
(whence \ 
\m{\Ext^1_{\ko_2}(\mm,\mm\ot_{\ko_1}\ki)\simeq\ker(\Psi_2)/\imm(\Psi_1)}) 
where \m{\Psi_1}, \m{\Psi_2} are represented by the matrices
\[\begin{pmatrix}y_0 & -x_0 \\ 0 & 0\\ 0 & 0\end{pmatrix} \ , \quad
\begin{pmatrix}0 & y_0 & -x_0 \\ 0 & 0 & 0\\ 0 & 0 & 0\end{pmatrix}
\]
respectively. It follows easily that we have
\[\ker(\Psi_2) \ = \ 
\{(\alpha,x_0\ot\ov{\epsilon},y_0\ot\ov{\epsilon}); \ \alpha\in \mm, \ 
\ov{\epsilon}\in\ki\} \ = \ (\mm\ot_{\ko_1}\ki)\oplus W \ , \]
with \ \m{W=(x_0,y_0)\ot\ki\subset 2(\mm\ot_{\ko_1}\ki)},
\[\imm(\Psi_1) \ = \ \mm^2\ot_{\ko_1}\ki \ \subset \ \mm\ot_{\ko_1}\ki \quad 
(\text{the 
first factor in} \ 3(\mm\ot_{\ko_1}\ki) ) \ . \]
Hence \ 
\m{\Ext^1_{\ko_2}(\mm,\mm\ot\ki_X)\simeq\big((\mm/\mm^2)\ot_{\ko_1}\ki\big)
\oplus\ko_1} \ .

\sepsubsub

{\bf Step 3. }{\em Description of the extensions -- } The morphisms \ 
\m{f:3\ko_2\to \mm\ot_{\ko_1}\ki} \ such that \Nligne \m{f(\imm(\phi_2))=\nsp} 
(i.e. the elements of \m{\ker(\Psi_2)}) are the \m{f_{\ov{\tau},\rho}}, with
\m{\ov{\tau}\in\mm\ot_{\ko_1}\ki}, \m{\rho\in\ko_1} and
\[f_{\ov{\tau},\rho}(a,b,c) \ = \ a_0\ov{\tau}+\rho(x_0b_0\ov{t}+y_0c_0\ov{t})
 \ , \]
inducing \ \m{\Theta_{\ov{\tau},\rho}:K\to\mm\ot_{\ko_1}\ki}. Let \ \m{\tau\in
(\mm/\mm^2)\ot_{\ko_1}\ki} \ be the image of $\ov{\tau}$. Then the element of 
\m{\Ext^1_{\ko_2}(\mm,\mm\ot\ki_X)} corresponding to \m{f_{\ov{\tau},\rho}} is 
\m{(\tau,\rho)}.

Let
\begin{equation}\label{equ4}0\lra\mm\ot_{\ko_1}\ki\hfl{i}{}M\hfl{p}{}\mm\lra 0
\end{equation}
be the extension associated to the image of \m{f_{\ov{\tau},\rho}} in 
\m{\Ext^1_{\ko_2}(\mm,\mm\ot_{\ko_1}\ki)}. We have an exact sequence
\begin{equation}\label{equ5} \xymatrix{0\ar[r] &  
K\ar[rrr]^-{g=(\iota,\Theta_{\ov{\tau},\rho})} & & &
2\ko_2\oplus(\mm\ot_{\ko_1}\ki)\ar[r]^-q & M\ar[r] & 0 \ ,}
\end{equation}
where $\iota$ is the inclusion. Hence we have
\begin{equation}\label{equ21}\imm(g) \ = \ 
\{(ey+\ov{\gamma},-ex+\ov{\delta},e_0\ov{\tau}+\rho(x_0\ov{\gamma}+
y_0\ov{\delta})); \ e\in\ko_2,\ \ov{\gamma},\ov{\delta}\in\ki\} \ . 
\end{equation}
In the exact sequence \ref{equ4}, \m{i:\mm\ot_{\ko_1}\ki\to M} \ comes from the 
inclusion \Nligne \m{\mm\ot_{\ko_1}\ki\hookrightarrow 
2\ko_2\oplus(\mm\ot_{\ko_1}\ki)}, and \ \m{q:M\to\mm} \ from \m{\phi_0}.

It follows easily from $(\ref{equ5})$ that $M$ is a \m{\ko_1}-module (i.e. 
\m{\ki M=\nsp}) if and only if \m{\rho=0}.

\sepsubsub

{\bf Step 4. }{\em Properties of the extensions -- } We will describe the 
morphism \ \m{h:\mm\ot_{\ko_1}\ki\to\mm\ot_{\ko_1}\ki} \ induced by the 
canonical morphism \ \m{M\ot_{\ko_2}\ki\to M} \ and $(\ref{equ4})$ ($M$ is 
balanced if and only if $h$ is an isomorphism). Let 
\m{\alpha\in\mm}, that we can write as \ \m{\alpha=\phi_0(a,b)=a_0x_0+b_0y_0}, 
with \m{a,b\in\ko_2}. Then from step 3, \m{h(\alpha\ot\ov{t})} is the unique 
\ \m{\beta\ot\ov{t}\in\mm\ot_{\ko_1}\ki} \ such that \ \m{(0,0,\beta\ot\ov{t})-
(\ov{t}a,\ov{t}b,0)\in g(K)}. Using $(\ref{equ5})$, we see that \ 
\m{\beta=-\rho\alpha}, so that $h$ is the multiplication by \m{-\rho}.

Now suppose that \m{\rho\not=0}. Then is is easy to see that \ 
\m{M^{(1)}=\imm(i)=\mm\ot_{\ko_1}\ki}, and that \ 
\m{M^{(1)}/M_1\simeq\mm/\rho\mm}. It follows that $M$ is balanced if and only 
if $\rho$ is invertible, and that $(\ref{equ4})$ is the canonical exact 
sequence \ \m{0\to M^{(1)}\to M\to M^{(0)}\to 0} \ if and only id \m{\rho=-1}. 

{\bf Step 5. }{\em Description of the morphisms \m{M\to\ko_2} -- } We suppose 
that \m{\rho=-1}.  Such a morphism comes from \ 
\m{\psi:2\ko_2\oplus\mm\ki\to\ko_2} \ such that \ \m{\psi\circ g=0}. There 
exist \m{a,b\in\ko_2}, \m{\alpha\in\ko_1}, such that \ 
\m{\psi(u,v,\ov{w})=au+bv+\alpha\ov{w}} \ for every \m{u,v\in\ko_2}, 
\m{\ov{w}\in\mm\ki}. Using $(\ref{equ21})$, we see that \ \m{\psi\circ g=0} \ 
is equivalent to
\[ax-by+\alpha\ov{\tau} \ = \ 0 \ , \quad a_0=\alpha x_0 \quad\text{and}\quad
b_0=\alpha y_0 \ . \]
Suppose that \m{\alpha=1}. We have then 
\[ay-bx \ = \ -\ov{\tau} \ , a-0 \ = \ x_0 \ , b_0 \ = y_0 \ . \]
So we can write \ \m{a=x+\ov{A}}, \m{b=y+\ov{B}}, with \ 
\m{\ov{A},\ov{B}\in\ki}, and \ \m{\ov{\tau}=-y_0\ov{A}+x_0\ov{B}}.

Now we show that $\psi$ is injective. Let \ \m{(u,v,\lambda)\in 2\ko_2\oplus
\mm\ki} \ be such that \ \m{\psi(u,v,\lambda)=0}, i.e. 
\m{u(x+\ov{A})+v(y+\ov{B})+\lambda=0}. Then we have \ \m{u_0x_0+v_0y_0=0}, 
hence we can write \ \m{u=ey+\ov{\gamma}}, \m{v=-ex+\ov{\delta}}, with 
\m{e\in\ko_2}, \m{\gamma,\delta\in\ki}, and we obtain \ 
\m{-e_0\ov{\tau}+x\ov{\gamma}+y\ov{\delta}+\lambda=0}, whence \ 
\m{(u,v,\lambda)=g(ey+\ov{\gamma},-ex+\ov{\delta})}, and $\psi$ is injective.

Since $\psi$ is injective, $M$ is isomorphic to \m{\imm(\psi)\subset\ko_2}, and 
\m{\imm(\psi)=(x+\ov{A},y+\ov{B})}. This proves theorem \ref{theo2}.
\end{proof}

\end{sub}

\sepprop

\Ssect{Parameterization of $\J_P$}{param_J}

Let \m{x,y\in\ko_2} be such that \m{x_0,y_0} are generators of $\mm$. If 
\m{\ov{A},\ov{B}\in\ki}, let \m{\boldsymbol{\tau}_{x,y}(\ov{A},\ov{B})} be the 
image of \ \m{-y_0\ov{A}+x_0\ov{B}} \ in \m{\mm\ki/\mm^2\ki}. By proposition 
\ref{prop23},
\[\xymatrix@R=5pt{\J_P\ar[r] & \mm\ki/\mm^2\ki\\ (x+\ov{A},y+\ov{B})\fmaps[r] 
& 
\boldsymbol{\tau}_{x,y}(\ov{A},\ov{B})}\]
is a bijection.

Let \m{x',y'\in\ko_2} be such that \m{x'_0,y'_0} are generators of $\mm$. Then 
there exist \ \m{\ov{u},\ov{v}\in\mm\ki} \ such that \ 
\m{(x',y')=(x+\ov{u},y+\ov{v})}. There exists an invertible matrix 
\m{\begin{pmatrix}\alpha & \beta\\ \gamma & \delta
\end{pmatrix}} with coefficients in \m{\ko_2}, such that
\[\begin{pmatrix}x'\\y'\end{pmatrix} \ = \ \begin{pmatrix}\alpha & \beta\\ 
\gamma & \delta\end{pmatrix}\begin{pmatrix}x+\ov{u}\\y+\ov{v}\end{pmatrix} \ .
\]
Let \ \m{\Delta=\alpha\delta-\beta\gamma}. If \m{\ov{A'},\ov{B'}\in\ki}, then \ 
\m{(x'+\ov{A'},y+\ov{B'})=(x+\ov{A},y+\ov{B})} \ with
\[\ov{A} \ = \ \frac{\delta_0}{\Delta_0}\ov{A'}-\frac{\beta_0}{\Delta_0}\ov{B'}
+\ov{u} \ , \quad \ov{B} \ = -\frac{\gamma_0}{\Delta_0}\ov{A'}+
\frac{\alpha_0}{\Delta_0}\ov{B'}+\ov{v} \ . \]
It follows that
\[\boldsymbol{\tau}_{x,y}(\ov{A},\ov{B}) \ = \ \frac{1}{\Delta_0}
\boldsymbol{\tau}_{x',y'}(\ov{A'},\ov{B'})-y_0\ov{u}+x_0\ov{v} \ , \]
We have \ \m{\mm\ki/\mm^2\ki\simeq\Omega_{X,P}\ot L_P}. Let \m{\Delta(x,y)} be 
the image of \m{x_0\wedge y_0} in \m{\omega_{X,P}}, and denote by 
\m{\Delta(x,y)^{-1}} the corresponding element of \m{\omega_{X,P}^*}. For 
\m{\ov{A},\ov{B}\in\ki} and \m{J=(x+\ov{A},y+\ov{B})}, let 
\[\boldsymbol{\lambda}_{P,x,y}(J) \ = \ 
\Delta(x,y)^{-1}\ot\boldsymbol{\tau}_{x,y}(\ov{A},\ov{B}) \ \in \ 
\omega_{X,P}^*\ot\Omega_{X,P}\ot L_P \ \simeq \ T_{X,P}\ot L \ . \]
Then \ \m{\boldsymbol{\lambda}_{P,x,y}:\J_P\to T_{X,P}\ot L} \ is a bijection.
We have
\begin{equation}\label{equ22}
\boldsymbol{\lambda}_{P,x',y'}(J) \ = \ \boldsymbol{\lambda}_{P,x,y}(J)+
\Delta(x,y)^{-1}\ot(-y_0\ov{u}+x_0\ov{v}) \ .
\end{equation}
Hence for every \m{J_1,J_2\in\J} we have
\[\boldsymbol{\lambda}_{P,x,y}(J_1)-\boldsymbol{\lambda}_{P,x,y}(J_2) \ = 
\boldsymbol{\lambda}_{P,x',y'}(J_1)-\boldsymbol{\lambda}_{P,x',y'}(J_2) \ . \]

It follows that

\sepprop

\begin{subsub}\label{prop24}{\bf Proposition: } $\J_P$ has a natural 
structure of affine space, with underlying $\C$-vector space \m{T_{X,P}\ot L}.
\end{subsub}

\end{sub}

\sepsub

\Ssect{Ideal sheaves on \m{X_2}}{Id_X2}

Let \m{J\in\J_P}. Then $J$ extends naturally to an ideal sheaf \m{J_2} on 
\m{X_2}, $P$ being the unique point \m{Q\in X} such that \ 
\m{J_{2,Q}\not=\ko_{X_2,Q}}. Hence we can see \m{\J_P} as a set of ideal 
sheaves on \m{X_2}.

We will study \ \m{\I=\dsp\bigcup_{P\in X}\kj_P}. We will see in 
\ref{uni_s} that $\I$ has a natural structure of an affine bundle on $X$, with 
associated vector bundle \m{T_X\ot L}. To describe this structure we will see 
here how it is related to the construction of \m{X_2} given in \ref{PMS}.

\sepprop

\begin{subsub}\label{ext_i} Extension of ideals -- \rm
Let $Z$ be a smooth variety, $U$ 
an open subset of \m{X_2}, and \m{U_0=U\cap X}. Let \ \m{\alpha:Z\to U_0} \ be 
a morphism. Then \ \m{Z_\alpha=\{(z,\alpha(z)); \ z\in Z\}} \ is a closed 
subvariety of \m{Z\times X}. Let $\ki$ be a sheaf of ideals on \m{Z\times U} 
such that the support of \m{\ko_{Z\times U}/\ki} (i.e. the set of points where 
\m{\ki\not=\ko_{Z\times U}}) is \m{Z_\alpha}. Then $\ki$ extends naturally to a 
sheaf of ideals \m{\ki'} on \m{Z\times X_2}, such that the support of 
\m{\ko_{Z\times X_2}/\ki'} is \m{Z_\alpha}.
\end{subsub}

\sepprop

\begin{subsub}\label{fam_id} Families of ideals -- \rm Suppose that
\begin{enumerate}
\item[--] The open subset \m{U_0} is affine.
\item[--] There exists a section \ \m{\sigma:\ko_X(U_0)\to\ko_{X_2}(U)} \ of 
the restriction \ \m{\ko_{X_2}(U)\to\ko_X(U_0)}, so we can view \m{\ko_X(U_0)} 
as a subring of \m{\ko_{X_2}(U)}.
\item[--] there exist \m{x,y\in\ko_X(U_0)} such that for every \m{Q\in U_0}, 
\m{\mm_Q=(x-x(Q),y-y(Q))} (where \m{\mm_Q} is the maximal ideal of 
\m{\ko_{X,Q}}), hence \m{\Omega_X(U_0)} is generated by \m{dx} and 
\m{dy}.
\item[--] $L_{|U}=\ki_{X|U}$ is trivial. Let $\ov{t}$ be a generator of 
\m{L_{|U}}.
\end{enumerate}
So we have a local description of \m{\Omega_{X|U_0}}:
\begin{enumerate}
\item[--] with differentials: \m{\Omega_X(U_0)} is generated by \m{dx} and 
\m{dy}.
\item[--] for every \m{Q\in U_0}, \m{\Omega_{X,Q}\simeq\mm_Q/\mm_Q^2}, and for 
every \m{\alpha,\beta\in\ko_{X,Q}}, \m{\alpha dx+\beta dy} corresponds to the 
image of \ \m{\alpha(Q)(x_0-x_0(Q))+\beta(Q)(y_0-y_0(Q))} \ in 
\m{\mm_Q/\mm_Q^2}. 
\end{enumerate}

Let \ \m{Y=U\times(U_0\times\C^2)}, \m{p:Y\to U}, \m{p_0:Y\to U_0},
\m{q_1,q_2:Y\to\C} be the projections. Then \ 
\m{(p^*(\sigma(x))-p_0^*(x_0)+q_1\ov{t},p^*(\sigma(y))-p_0^*(y_0)+q_2\ov{t})} \ 
is a regular 
sequence at each point of $Y$. Let \ 
\[{\bf I} \ = \ 
(p^*(\sigma(x))-p_0^*(x_0)+q_1\ov{t},p^*(\sigma(y))-p_0^*(y_0)+q_2\ov{t}) \ , 
\]
which is an ideal sheaf on $Y$, of a subscheme \m{Z\subset Y}. One of the 
irreducible and connected components of $Z$, \m{Z_0}, contains the points 
\m{(P,P,u)}, \m{P\in U_0,u\in\C^2}. Let \ 
\m{\boldsymbol{\ki}^{\sigma,\ov{t},x,y}=\ki_{Z_0}}, which by \ref{ext_i} can be 
seen as a family of ideal sheaves on \m{X_2} parameterized by \m{U_0\times\C_2} 
and flat on it. For every \ \m{(P,a,b)\in U_0\times\C^2}, we have
\[\big(\boldsymbol{\ki}^{\sigma,\ov{t},x,y}\big)_{P,a,b} \ = \ 
(\sigma(x)-x(P)+a\ov{t},\sigma(y)-y(P)+b\ov{t}) \ . \]

We have an isomorphism
\[\Theta^{x,y,\ov{t}}:\ko_{U_0}\ot\C^2\lra
\big(\Omega_{X}\ot\omega_{X}^*\ot L\big)_{|U_0} \ = \ (T_X\ot L)_{|U_0} \ , \]
defined by
\[\Theta^{x,y,\ov{t}}(a,b) \ = \ (-a.dy+b.dx)\ot(dx\wedge dy)^{-1}\ot\ov{t} \ 
, \]
for every \m{a,b\in\ko_X(U_0)}. Let
\[\boldsymbol{\kj}^{\sigma,\ov{t},x,y} \ = \ 
((\Theta^{x,y,\ov{t}})^{-1})^\sharp(\boldsymbol{\ki}^{\sigma,\ov{t},x,y}) \ , 
\]
so that \ 
\m{\big(\boldsymbol{\kj}^{\sigma,\ov{t},x,y}\big)_{P,\Theta^{x,y,\ov{t}}(a,b)}=
\big(\boldsymbol{\ki}^{\sigma,\ov{t},x,y}\big)_{P,a,b}} .

\sepprop

\begin{subsub}\label{prop32}{\bf Proposition: } 
\m{\boldsymbol{\kj}^{\sigma,\ov{t},x,y}} \ is independent of $x$,$y$ and 
$\ov{t}$.
\end{subsub}
\begin{proof}
For every \m{Q\in U_0} we have
\[\big(\Theta^{x,y,\ov{t}}(a,b)\big)_Q \ = \ 
\boldsymbol{\lambda}_{Q,\sigma(x),\sigma(y)}
\big((\boldsymbol{\kj}^{\sigma,\ov{t},x,y})_{Q,\Theta^{x,y,\ov{t}}(a,b)}\big) 
\ . \]
In other words, for every \ \m{\eta\in\big(\Omega_{X}\ot\omega_{X}^*\ot 
L\big)_Q}, 
\[\boldsymbol{\lambda}_{Q,\sigma(x),\sigma(y)}
\big((\boldsymbol{\kj}^{\sigma,\ov{t},x,y})_{Q,\eta}\big) \ = \ \eta \ . \]
If \m{x,y} are replaced with \m{x',y'\in\ko_X(U_0)}, we have, for every \m{P\in 
U_0}
\[(\sigma(x)-x(P),\sigma(y)-y(P)) \ = \ (\sigma(x')-x'(P),\sigma(y')-y'(P))\]
(it is the ideal of \m{\ko_{X_2,P}} generated by \m{\sigma(\mm_Q)}). Hence, by 
formula $(\ref{equ22})$, the maps 
\m{\boldsymbol{\lambda}_{Q,\sigma(x),\sigma(y)}} and 
\m{\boldsymbol{\lambda}_{Q,\sigma(x'),\sigma(y')}} are the same, and since they 
are isomorphisms, \m{\boldsymbol{\kj}^{\sigma,\ov{t},x,y}} is independent of 
the choice of $x$, $y$ and $\ov{t}$. 
\end{proof}

\sepprop

We will note \ 
\m{\boldsymbol{\kj}^\sigma=\boldsymbol{\kj}^{\sigma,\ov{t},x,y}}. It is a flat 
family of ideal sheaves on \m{X_2}, parameterized by \m{(T_X\ot L)_{|U_0}}.

The map
\[\xymatrix@R=5pt{f_Q:T_{X,Q}\ot L_Q=\big(\Omega_{X}\ot\omega_{X}^*\ot 
L\big)_Q\ar[r] & \J_Q\\ \eta\fmaps[r] & (\boldsymbol{\kj}^\sigma)_{Q,\eta}}\]
is a bijection, i.e. \m{\big(\Omega_{X}\ot\omega_{X}^*\ot 
L\big)_Q} \ parameterizes the set of ideal sheaves of \m{\ko_{X_2}}, which are 
balanced and whose restriction to $X$ is the ideal sheaf of $Q$.

Let \ \m{\J_{U_0}=\bigcup_{Q\in U_0}\J_Q}. Then \m{\J_{U_0}} parameterizes the 
set of ideal sheaves of \m{\ko_{X_2}}, which are balanced and whose restriction 
to $X$ is the ideal sheaf of a point of \m{U_0}. We have a bijection
\[\mu_\sigma: (T_X\ot L)_{|U_0}\lra\J_{U_0}\]
over \m{U_0} defined by the maps \m{f_Q}.
\end{subsub}

\sepsubsub

\begin{subsub}\label{c_par} -- Change of parameters -- \rm
Suppose that we have another section \ \m{\sigma':\ko_X(U_0)\to\ko_{X_2}(U)} \ 
of the restriction \m{\ko_{X_2}(U)\to\ko_X(U_0)}. Then there is a derivation 
$D$ of \m{\ko_X(U_0)}such that
\[\sigma'(\lambda) \ = \ \sigma(\lambda)+D(\lambda)\ov{t} \ , \]
for every \m{\lambda\in\ko_X(U_0)}.
Let
\[x' \ = \ \sigma'(x) \ = \ \sigma(x)+D(x)\ov{t} \ , \qquad y' \ = \ \sigma'(y) 
\ = \ \sigma(y)+D(y)\ov{t} \ . \]
We can view \m{D\ot\ov{t}} as a section of \m{(T_X\ot L)_{|U_0}}.

\sepprop

\begin{subsub}\label{prop25}{\bf Proposition: } For every \m{Q\in U_0}, \m{u\in 
T_{X,Q}\ot L_Q},
\[\mu_{\sigma'}(u) \ = \ \mu_\sigma\big(u+(D\ot\ov{t})(Q)\big) 
\ . \]
\end{subsub}
\begin{proof}
From the definitions, using \m{x'}, \m{y'} instead of $x$, $y$, we have
\begin{eqnarray*}
\big(\boldsymbol{\kj}^{\sigma'}\big)_{P,\Theta^{x',y',\ov{t}}(a,b)} & = & 
(\sigma(x)-x(P)+(a+D(P)(x))\ov{t},\sigma(y)-y(P)+(b+D(P)(y))\ov{t})  \ ,\\
& = & 
\big(\boldsymbol{\kj}^\sigma\big)_{P,\Theta^{x,y,\ov{t}}(a+D(P)(x),b+D(P)(y))} 
\ , \end{eqnarray*}
which is equivalent to
\[\mu_{\sigma',\ov{t}}((-a.dy+b.dx)\ot(dx\wedge dy)^{-1}\ot\ov{t}) \ = \ 
\qquad\qquad\qquad\qquad\qquad\qquad\]
\[\qquad\qquad\qquad\qquad
\mu_{\sigma,\ov{t}}\big(-(a+D(x))dy+(b+D(y))dx\big)\ot(dx\wedge 
dy)^{-1}\ot\ov{t} \ . \]
Using the canonical isomorphism \ \m{\omega_{X,P}^*\ot\Omega_{X,P}\simeq 
T_{X,P}}, one obtains proposition \ref{prop25}.
\end{proof}

\end{subsub}

\sepprop

\end{sub}

\sepsub

\Ssect{Moduli spaces of ideals and universal sheaves}{uni_s}

Using the results of \ref{Id_X2} we show that the set of ideals $\I$ has a 
natural structure of an affine bundle on $X$, with associated vector bundle 
\m{T_X\ot L}. We construct a {\em universal sheaf} on \ \m{\I\times X_2}.

\sepprop

\begin{subsub}\label{idx} Ideal sheaves on $X$ -- \rm Let \ \m{\Delta\subset 
X\times X} be the diagonal. Then \m{\ki_\Delta} can be seen as a family of 
ideal sheaves on $X$, parameterized by $X$. For every \m{x\in X}
\[\ki_{\Delta,x} \ = \ \ki_{\Delta|X\times\{x\}} \ = \ \ki_x \ . \]
\end{subsub}

\sepprop

\begin{subsub}\label{mod_s} Construction of the moduli space and the universal 
sheaf -- \rm
Let \m{(V_i)_{1\leq i\leq k}} be an open affine cover of $X$ such that for each 
\m{V_i} the conditions of \ref{fam_id} are satisfied: we have a section \ 
\m{\sigma_i:\ko_X(V_i)\to\ko_{X_2}(V_i)} \ of the restriction \ 
\m{\ko_{X_2}(V_i)\to\ko_X(V_i)}, and \m{\ov{t}_i} is a generator of 
\m{L_{|V_i}}. We obtain the family \m{\kj^{\sigma_i}} of ideal sheaves on 
\m{X_2} parameterized by \m{(T_X\ot L)_{|V_i}}. If \ \m{1\leq i<j\leq k} \ we 
have
\[\sigma_{j|V_{ij}} \ = \ \sigma_{i|V_{ij}}+D_{ij}.\ov{t}_i \ , \]
where \m{D_{ij}} is a derivation of \m{\ko_X(V_{ij})}. With the notations of 
\ref{PMS}, \m{(\alpha_i^{-1}\ot D_{ij})_{1\leq i<j\leq k}} 
represents an element $\lambda$ of \m{H^1(X,T_X\ot L)}. From \ref{DBC},
\m{\P(H^1(X,T_X\ot L))\cup\{0\}} \ parameterizes the set of primitive double 
schemes $Y$ such that \m{Y_{red}=X} and associated line bundle $L$, and \ 
\m{\C\lambda=\zeta(X_2)} is the element corresponding to \m{X_2}.

We have an affine isomorphism
\[\xymatrix@R=5pt{\rho_{ij}:(T_X\ot L)_{|V_{ij}}\ar[r] & (T_X\ot L)_{|V_{ij}}\\
u\fmaps[r] & u+D_{ij}(u)\ot\ov{t}_i}\]
and
\[\rho_{ij}^\sharp(\boldsymbol{\kj}^{\sigma_i}) \ \simeq \ 
\boldsymbol{\kj}^{\sigma_j} \ . \]
The \m{\rho_{ij}} define an affine bundle $\I$ on $X$, with associated vector 
bundle \m{T_X\ot L}. According to \ref{aff_b} we have

\sepprop

\begin{subsub}\label{prop35}{\bf Proposition: } The element \m{\eta(\I)} of \ 
\m{\big(\P(H^1(X,T_X\ot L))\cup\{0\}\big)/\Aut(T_X\ot L)} \ associated to $\I$ 
is the image of \m{\zeta(X_2)}.
\end{subsub}

\sepprop

The closed points of $\I$ are the ideal sheaves $\ki$ on \m{X_2} of subschemes 
$Z$ such that
\begin{enumerate}
\item[--] $Z$ contains only one closed point $P$ of $X$.
\item[--] there exist \m{x,y\in\ko_{X_2,P}} such that their images in 
\m{\ko_{X,P}} generate the maximal ideal, and \ \m{\ki_P=(x,y)}.
\end{enumerate}

The sheaves \m{\boldsymbol{\kj}^{\sigma_i}} can be glued, using \ref{c_par}, to 
obtain a flat family of ideal sheaves $\boldsymbol{\kj}$ parameterized by $\I$. 
It is a balanced sheaf on \m{X_2\times\I}. If \ \m{p_X:\I\to X} \ is the 
projection, we have
\[\boldsymbol{\kj}_{|X\times\I} \ \simeq \ p_X^\sharp(\ki_\Delta) \ . \]
\end{subsub}

\end{sub}

\sepsub

\Ssect{Moduli spaces of sheaves}{shea}

The ideal sheaves \m{\ki\in\I} can also be deformed by tensoring with a line 
bundle. This is why one must consider also the 
sheaves \m{\ki\ot D}, \m{\ki\in\I}, \m{D\in\Pic(X_2)}. In this way we obtain 
here a moduli space of sheaves (with a suitable universal property described in 
proposition \ref{prop33}), isomorphic to \ \m{\I\times\Pic(X_2)}.

\sepprop

\begin{subsub}\label{mod_X} Moduli spaces on rank 1 sheaves on $X$ -- \rm Let 
$\D$ be a Poincar\'e bundle on \m{\Pic(X)} and \ \m{\M_X=X\times\Pic(X)}. 
Let \m{p_1:\M_X\to X}, \m{p_2:\M_X\to\Pic(X)} be the projections, and
\[\mathfrak{D} \ = \ p_1^\sharp(\ki_\Delta)\ot p_2^\sharp(\D) \ . \]
It is a flat family of sheaves on $X$ parameterized by \m{\M_X}. For every \ 
\m{(x,D)\in\M_X}, \m{\mathfrak{D}_{(x,D)}=\ki_x\ot D}.
\end{subsub}
The following properties are well known:
\begin{enumerate}
\item[--] Let \ \m{(x,D),(x',D')\in\M_X}. Suppose that \ 
\m{\mathfrak{D}_{(x,D)}\simeq\mathfrak{D}_{(x',D')}}. Then \ \m{(x,D)=(x',D')}. 
So we can see \m{\M_X} as the family of sheaves \m{\mathfrak{D}_{(x,D)}}.
\item[--] Let $\kf$ be a flat family of sheaves on $X$ parameterized by a 
scheme 
$Y$. Let \m{y\in Y} be a closed point such that \ \m{\kf_y\in\M_X}. Then 
there exists an open neighborhood $V$ of $y$ such that for every closed point 
\m{v\in V}, \m{\kf_v\in\M_X}.
\item[--] Let $\kf$ be a flat family of sheaves on $X$ parameterized by a 
scheme 
$Y$. Suppose that for every closed point \m{y\in Y}, \m{\kf_y\in\M_X}. Then 
there exists a unique morphism \ \m{f_\kf:Y\to\M_X}, and a line bundle 
\m{\Gamma} on $Y$ such that \ \m{\kf\simeq f^\sharp(\mathfrak{D})\ot 
p_Y^*(\Gamma)} (where \m{p_Y=X\times Y\to Y} is the projection).
\end{enumerate}
Thus \m{\M_X} is a moduli space of rank 1 sheaves on $X$ (cf. \ref{fin-mod}).

\sepprop

\begin{subsub}\label{mod_X2} Moduli spaces on rank 1 sheaves on \m{X_2} -- \rm 
We can do the same work on \m{X_2}, using $\I$ and \m{\Pic(X_2)}.
\end{subsub}

\begin{subsub}\label{lem21}{\bf Lemma:} Let $\ki$, \m{\ki'} be ideals in $\I$, 
and $\D$ a line bundle on \m{X_2}. If \ \m{\ki'\simeq\ki\ot\D} \ then 
\m{\D=\ko_{X_2}} and \ \m{\ki'=\ki}.
\end{subsub}
\begin{proof} This follows easily from the fact that \ 
\m{\ki^*={\ki'}^*=\ko_{X_2}}.
\end{proof}

\sepprop

Let \ \m{\boldsymbol{I}(X_2)=\I\times\Pic(X_2)}. Let 
\m{q_1:\boldsymbol{I}(X_2)\to\I}, \m{q_2:\boldsymbol{I}(X_2)\to\Pic(X_2)} \ be 
the projections. Let \m{(P_i)_{i\in I}} be an open cover of \m{\Pic(X_2)} such 
that for every \m{i\in I}, there is a Poincaré bundle \m{\kd_i} on \m{X_2\times 
P_i}, which define \m{\Pic(X_2)} as a {\em fine moduli space} of sheaves (cf. 
\ref{fin-mod}, \ref{pic-sch}).

Let \ \m{\boldsymbol{\kh_i}=q_1^\sharp(\boldsymbol{\kj})\ot 
q_2^\sharp(\kd_i)}, which is a sheaf on \ \m{\I\times P_i}. For  every ideal 
$\ki$ in $\I$ and \ \m{D\in P_i}, \m{\boldsymbol{\kh}_{\ki,D}=\ki\ot\kd_{i,D}} .

\sepprop

\begin{subsub}\label{prop33}{\bf Proposition: } Let $Y$ be a scheme and $\ke$ a 
coherent sheaf on \m{X_2\times Y}, flat over $Y$. Let \m{y\in Y} be a closed 
point such that \ \m{\ke_y\in\boldsymbol{I}(X_2)}. Then there exists an 
neighborhood \m{V\subset Y} of $y$ such that
\begin{enumerate}
\item[--] $\ke_{1|X\times V}$ and \m{\ke_{|X\times V}} are flat on $V$.
\item[--] For every closed point \m{v\in V}, \m{\ke_v\in\boldsymbol{I}(X_2)}, 
and \m{(\ke_1)_v=(\ke_v)_1}, \m{(\ke_{|X\times V})_v=(\ke_v)_{|X}}.
\end{enumerate}
\end{subsub}
\begin{proof} This follows from theorem \ref{theo5}, 1-.
\end{proof}

\sepprop

The proof of the following result is similar of that of theorem 1.1.3 of 
\cite{dr11}, and shows that \m{\boldsymbol{I}(X_2)} is a moduli space of 
sheaves in the sense of \ref{fin-mod}.

\sepprop

\begin{subsub}\label{theo6}{\bf Theorem: } Let $Y$ be a scheme and $\ke$ a 
coherent sheaf on \m{X_2\times Y}, flat over $Y$, such that for every closed 
point \m{y\in Y}, \m{\ke_y\in\boldsymbol{I}(X_2)}. Then there exists a unique 
morphism \ \m{f_\ke:Y\to \boldsymbol{I}(X_2)} \ and an open cover 
\m{(Y_i)_{i\in I}} such that for every \m{i\in I}, 
\m{f_\ke^\sharp(\mathfrak{\kh})_{|X_2\times Y_i}\simeq\ke_{|X_2\times Y_i}}.
\end{subsub}

\sepprop

Let \m{\Gamma(X_2)\subset\Pic(X)} be the algebraic subgroup of vector bundles 
that can be extended to \m{X_2}. According to \ref{pic-sch}, \m{\Pic(X_2)} is 
an affine bundle on \m{\Gamma(X_2)} with associated vector bundle \ 
\m{\ko_{\Gamma(X_2)}\ot H^1(X,L)}. Let \ \m{p_X:X\times\Gamma(X_2)\to X}, 
 \ be the projection. Then 
\m{\boldsymbol{I}(X_2)} is an affine bundle on \m{X\times\Gamma(X_2)}, with 
associated vector bundle \ \m{(\ko_{X\times\Gamma(X_2)}\ot H^1(X,L))\oplus
p_X^*(T_X\ot L)}.

\end{sub}

\vskip 4cm

\end{document}